\documentclass[11pt]{article}
\usepackage{epsfig,amssymb,amsmath}
\usepackage{amsfonts}
\usepackage{graphicx}
\usepackage{url}

\newcommand{\bif}{\; , \;}
\newcommand{\bis}{\; ; \;}

\newcommand{\sta}{{\rm sta}}

\newcommand{ \PP}{{P}}
\newcommand{ \LL}{{L}}

\newcommand{\vsig}{\varsigma}

\newcommand{\ww}{w}

\newcommand{\uu}{u}

\newcommand{\real}{{\mathbb R}} 

\newcommand{\half}{\frac{1}{2}}

\newcommand{\bsig}{\mbox{\boldmath$\sigma$}}
\newcommand{\bvsig}{\mbox{\boldmath$\varsigma$}}

\newcommand{\bveps}{\mbox{\boldmath$\xi$}}

\newcommand{\barbvsig}{\bar{\mbox{\boldmath$\varsigma$}}}

\newcommand{\beps}{\mbox{\boldmath$\epsilon$}}

\newcommand{\bLam}{{\mbox{\boldmath$\Lambda$}}}

\newcommand{\bchi}{{\mbox{\boldmath$\chi$}}}

\newcommand{\barbchi}{\bar{\mbox{\boldmath$\chi$}}}
\newcommand{\bxi}{{\mbox{\boldmath$\xi$}}}

\newcommand{\barbxi}{\bar{\mbox{\boldmath$\xi$}}}

\newcommand{\la}{\langle}
\newcommand{\ra}{\rangle}

\newcommand{\bp}{{{\bf p}}}
\newcommand{\Oo}{\Omega}

\newcommand{\Gt}{{\Gamma_t}}
\newcommand{\Gu}{{\Gamma_u}}
\newcommand{\bE}{{\bf E}}
\newcommand{\bQ}{{\bf Q}}

\newcommand{\eba}{\begin{array}}
\newcommand{\eea}{\end{array}}
\newcommand{\ebe}{\begin{eqnarray}}
\newcommand{\eee}{\end{eqnarray}}
\newcommand{\eb}{\begin{equation}}
\newcommand{\ee}{\end{equation}}

\newcommand{\calW}{{\cal{W}}}
\newcommand{\calP}{{\cal{P}}}

\newcommand{\bg}{{\bf g}}
\newcommand{\bC}{{\bf C}}
\newcommand{\bG}{{\bf G}}
\newcommand{\bH}{{\bf H}}

\newcommand{\bt}{{\bf t}}
\newcommand{\bff}{{\bf f}}

\newcommand{\bR}{{\bf R}}
\newcommand{\bq}{{\bf q}}
\newcommand{\bv}{{\bf v}}
\newcommand{\bw}{{\bf w}}

\newcommand{\bx}{{\bf x}}
\newcommand{\by}{{\bf y}}

\newcommand{\bN}{{\bf N}}
\newcommand{\bF}{{\bf F}}
\newcommand{\bA}{{\bf A}}

\newcommand{\bI}{{\bf I}}

\newcommand{\calS}{{\cal S}}

\newcommand{\calL}{{\cal L}}

\newcommand{\calC}{{\cal C}}

\newcommand{\calZ}{{\cal Z}}
\newcommand{\calU}{{\cal U}}
\newcommand{\calE}{{\cal E}}

\newcommand{\calV}{{\cal V}}

\newcommand{\calY}{{\cal Y}}
\newcommand{\calX}{{\cal X}}

\newcommand{\bu}{{\bf u}}

\newcommand{\barbu}{\bar{\bf u}}

\newcommand{\barbx}{\bar{\bf x}}

\newcommand{\bs}{{\bf s}}

\newcommand{\dt}{{\,\mbox{d}t}}

\newcommand{\be}{{\bf e}}

\newcommand{\dO}{\,\mbox{d}\Oo}

\newcommand{\dG}{\,\mbox{d} \Gamma}

\newcommand{\alp}{{\alpha}}

\newcommand{ \sig}{{\sigma}}
\newcommand{ \Lam}{{\Lambda}}

\newcommand{ \lam}{{\lambda}}
\newcommand{ \xx}{{\bf x}}

\newcommand{\bD}{{\bf D}}

\newcommand{\tr}{{\mbox{tr}}}

\newcommand{\rank}{{\mbox{rank }}}
\newcommand{\Diag}{{\mbox{Diag }}}

\newtheorem{lemma}{Lemma}
\newtheorem{assumption}{Assumption}
\newtheorem{remark}{Remark}

\newtheorem{thm}{Theorem}

\newtheorem{definition}{Definition}

\renewcommand\uu{{\bf u}}

\newcommand\UU{{U}}
\newcommand\DD{{D}}
 \newcommand\WW{{W}}
 \newcommand\VV{{\Phi}}

    \newcommand\bz{{\bf z}}
  \renewcommand\ww{{\bf w}}
   \newcommand\TT{{T}}
 \newcommand\FF{{F}}
 \newcommand\G{G\^{a}teaux }
 \newcommand\z{Z\u{a}linescu}
  \newcommand\za{Z\u{a}linescu }
 \newcommand\Gap{G_{ap}}
 \newcommand\bGd{\bG_{\delta_k}}
 \newcommand\bfd{\bff_{\delta_k}}

\newtheorem{Conj}{Conjecture}
 \renewcommand{\bveps}{\mbox{\boldmath$\varepsilon$}}
\begin{document}
\begin{center}
 {\Large{\bf  On  Unified Modeling, Canonical Duality-Triality Theory,\\ Challenges and Breakthrough in Optimization 
 }} \vspace{.5cm}\\

{\bf David Yang Gao }\\
{\small {\em
  {\small  Faculty of Science and  Technology,\\
Federation University,
Mt Helen, Victoria 3353, Australia}
}}
\end{center}

\begin{abstract}
 A unified model  is addressed  for  general optimization problems in multi-scale complex systems.
 Based on   necessary conditions and basic principles in physics,
the canonical duality-triality theory is presented in a precise way to include
  traditional duality theories  and popular methods as special applications.
  Two conjectures on NP-hardness are discussed, which should play important roles for correctly understanding and
  efficiently solving challenging real-world problems.
  Applications are illustrated
  for  both nonconvex continuous optimization and mixed  integer nonlinear programming.
  Misunderstandings and confusion on some basic concepts, such as objectivity, nonlinearity,
 Lagrangian, and  Lagrange multiplier method are discussed and classified.
Breakthrough from recent false  challenges   by
  C.  Z\u{a}linescu and his co-workers  are addressed.
This paper will bridge a significant gap between optimization and
 multi-disciplinary fields of  applied math  and computational sciences.

\end{abstract}


\section{Introduction and Motivation}
 General problems in mathematical optimization  are  usually formulated in   the   following   form
\eb
\min f(\xx) , \;\; \mbox{ s.t. } \;  \bg (\xx) \le 0   , \label{eq-go}
\ee
where the unknown $\xx\in \real^n$ is a vector,
  $f(\xx):\real^n \rightarrow \real$ is the so-called ``objective" function\footnote{This terminology is used mainly
  in English literature.
The function $f(\xx)$ is called  the target function in Chinese and Japanese literature.}, %
and    $ \bg (\xx) = \{ g_j(\xx) \} : \real^n \rightarrow \real^m$  is a vector-valued  constraint function.
It must be emphasized  that,  different from the basic concept of {\em objectivity} in continuum physics and nonlinear analysis,
the objective function used extensively in  optimization literature
 is allowed to be any arbitrarily given function,  even the linear function.
Therefore,  this mathematical problem  is artificial. Although it enables   one to ``model" a very wide range of problems,
 it comes at a price: many  
 global  optimization problems   are considered to be NP-hard.
Without detailed  information on these arbitrarily given functions,
it is impossible to have a powerful  theory for solving  the general nonconvex problem (\ref{eq-go}).

Canonical duality-triality  is a newly developed and  continuously improved   methodological theory.
This theory comprises mainly  1) a canonical transformation, which is a  versatile methodology
 that can be used to model complex systems within a unified framework;
2) a complementary-dual principle, which can be used to formulate  a perfect dual problem with
 a unified analytic solution; and 3)  a triality theory, which can identify both global and local extrema and to develop effective canonical dual algorithms for solving
real-world  problems in both continuous and discrete systems.
This theory was developed from Gao and Strang's original work on nonconvex variational/boundary value problems in
large deformation mechanics \cite{gao-strang89}.
 It was shown in Gao's book \cite{gao-dual00} and
in  the recent articles \cite{lat-gao-opl,g-r-s}  that both the (external) penalty and
Lagrange    multiplier methods are special applications of
 the canonical duality theory in convex optimization.
 It is now understood that this theory
 reveals an intrinsic multi-scale duality pattern in complex systems,
 many popular theories and powerful methods in  nonconvex analysis, global optimization, and
 computational science can be unified within the framework of the  canonical duality-triality  theory.
Indeed, it is easy to show  that the  popular {\em semi-definite programming} (SDP) methods in global optimization
and the {\em half-quadratic regularization} in image processing  are naturally covered by  the canonical duality theory \cite{gao-bridge,lat-gao-half,z-g-y}.

Mathematics and mechanics have been complementary partners since the Newton times.
Many fundamental ideas, concepts, and mathematical methods extensively used in
calculus of variations and    optimization are originally developed from mechanics.
It is known that  the classical  Lagrangian duality theory and   the associated Lagrange multiplier method
 were  developed by Lagrange in analytical mechanics \cite{Lag}.
 The modern concepts of super-potential and sub-differential in convex analysis were proposed by J.J. Moreau from frictional
 mechanics \cite{more-68}.
However, as V.I. Arnold indicated  \cite{arnold}:
 ``In the middle of the twentieth century it was attempted to divide physics and mathematics. The consequences turned out to be catastrophic".

Indeed, the canonical duality theory was developed from some
 fundamental concepts   of objectivity and work-conjugate principle in continuum physics.
Due to the existing gap between nonlinear analysis/mechanics and optimization, this  theory has been
mistakenly challenged by  M. Voisei, C.  Z\u{a}linescu and his former student in a set of more than  12 papers.
Although knowledgeable scholars can easily understand the conceptual  mistakes in these challenges,
the large number of so-called ``counterexamples" and false conclusions
 generated negative impacts to the communities.
Therefore, it is necessary to have a formal response to the recent paper by \za \cite{z16}.
Instead of directly  responding all these challenges one by one, 
this paper presents    the canonical duality theory in a  systematical way from
a unified modeling,  basic  assumptions to the  theory, method, and applications.
The methodology,  conjectures, and responses
 are important  for understanding not only this  unconventional theory,
but also many  challenging problems in complex systems.
A hope to
 bridge the existing  gap between the mathematical optimization and the interdisciplinary fields of  mathematical physics is  author's main goal of this paper.

 \section{Multi-Scale Modeling and Properly Posed Problems} 




The  canonical duality theory
was developed from   Gao and Strang's   work on minimum potential energy principle for solving
 the following variational  problem 
  in large  deformation theory \cite{gao-strang89}:
\eb
(\calP_0): \;\;\; \min \{ \Pi (\uu) = \WW(D\uu) +  F(\uu) \; | \; \uu \in \calU_c   \}, \label{eq-gs}
\ee
where the unknown $\uu(\bx)$ is a 
function in a differentiable manifold
$\calU$  over a field; 
 $F:\calU_a \subset \calU\rightarrow \real$ is an {\em external energy};
  $D:\calU_a  \rightarrow \calW$ is a linear differential operator which assigns each configuration $\uu$
  to an internal variable
 $\bw = D \uu$ in different scale;
  the real-valued function
  $\WW:\calW_a\subset \calW  \rightarrow \real$ is  the so-called {\em internal (or stored) energy}.
  In $\calU_a$ the geometrical constraints (such as boundary and initial conditions) are pre-described for each given problem; while
  $\calW_a$ contains certain  physical (constitutive) constraints of the system.
The feasible set $\calU_c = \{ \uu \in \calU_a| \;\; D\uu \in \calW_a\} $ is known as the   {\em kinetically  admissible space}
in nonlinear field theory.

\subsection{ Objectivity, Subjectivity and Well-posed Problem}

 Objectivity is  a basic concept  in mathematical modeling and nonlinear analysis
 \cite{gao-dual00,ciarlet88,ciarlet,mar-hug}.
 Let  $ {\cal R} $  be
 a  special orthogonal group, i.e. $\bR \in {\cal R} $ if and only if
 $  \bR^T = \bR^{-1} $ and $ \det \bR = 1$.
 A  mathematical definition was given in Gao's book (Definition 6.1.2 \cite{gao-dual00}).
  \begin{definition}[Objectivity] 

 A set $\calW_a   $ is said to be objective if
 $
 \bR \bw \in \calW_a \;\; \forall \bw \in \calW_a, \; \forall \bR \in {\cal R}.
$
 A real-valued function $\WW:\calW_a \rightarrow \real$ is said to be objective if \vspace{-.1cm}
 \eb
 \WW(\bR \bw ) = \WW(\bw) \;\; \forall \bw \in \calW_a, \; \forall \bR \in {\cal R}.\vspace{-.1cm}
 \ee

 \end{definition}
\begin{lemma} A real-valued function $\WW(\bw)$ is objective if and only if there exists
a real-valued function $\Phi(\bE)$ such that $\WW(\bw) = \Phi(\bw^T\bw)$.\vspace{-.2cm}
\end{lemma}

 Geometrically speaking, an objective function is  rotational symmetry, which should be a SO$(n)$-invariant in $n$-dimensional Euclidean space.
  Physically, an objective function doesn't depend on observers.
   Because of Noether's theorem\footnote{i.e., every differentiable symmetry of the action of a physical system has a corresponding conservation law.}, rotational symmetry of a physical system is equivalent to the angular momentum conservation law (see Section 6.1.2 \cite{gao-dual00}).
   Therefore, the objectivity is essential for  any  real-world mathematical
 models.
It was  emphasized by P.  Ciarlet that the objectivity is not an assumption, but an axiom \cite{ciarlet88}.
  Indeed,  the objectivity is also known as the
  {\em axiom of  frame-invariance} \cite{ball}.

In Gao and Strang's work, 
 the internal energy  $W(\bw)$  must  be an  objective function   such that its variation (\G
     derivative)
    $\bsig = \partial  W(\bw)$ is  the so-called   {\em constitutive duality law},
which depends only on the intrinsic property of the system.
Dually, the external energy $F(\uu)$ can be  called the {\em subjective function} \cite{gao-ilk}, which depends on each  problem
    such that its variation is governed by the {\em action-reaction duality law}:
  $\barbu^* = - \partial \FF(\bu) \in \calU^*$.
A system is  conservative if the action  is independent of the reaction.
Therefore, the subjective function  must be linear on its domain $\calU_a$ and,
 by Riesz representation theorem,  we should have
  $F(\uu) = - \la \uu, \barbu^* \ra$, where
 the bilinear form $\la \uu, \uu^* \ra :\calU \times \calU^* \rightarrow \real$ puts $\calU$ and $\calU^*$ in duality.
  Together, $\Pi(\uu) = W(D\uu) + F(\uu)$ is called the total potential  and the minimum potential energy principle
  leads to the general variational problem (\ref{eq-gs}).
 From the point view of linguistics,  if we consider $\FF(\uu)$ as a subject,  $\WW(\ww)$ as an object,
 and the operation ``$+$" as a predicate,
  then $\Pi(\uu)$  forms   a grammatically correct  sentence.\footnote{By the facts that (object, subject) is a duality pair in a noun (or  pronoun)  space, 
  which is dual to a  verb space,  the multi-level duality pattern  $ \{ (object, subject) ; predicate\} $
  is called  triality,  which  is essential for  languages and sciences.}
  The criticality condition $\partial \Pi(\uu) = 0$ leads to the equilibrium equation
  \eb
 A (\uu) = \DD^* \partial \WW(\DD \uu) = \barbu^* \label{eq-equil}
  \ee
  where $\DD^*: \calW^*_a \rightarrow \calU^*$ is an adjoint operator of $\DD$ and $A:\calU_c \rightarrow \calU^*$ is
  called {\em equilibrium operator}.
The triality structure  $\mathbb{S}^e = \{ \la \calU, \calU^* \ra ;  A \} $
    forms  an elementary system  in Gao's book (Chapter 4.3, \cite{gao-dual00}).
  This abstract form   covers most well-known  equilibrium problems in  real-world applications ranging from
      mathematical physics in continuous analysis  to  mathematical programming in discrete  
      systems (see the celebrate text by Gil Strang \cite{strang}).
  Particularly, if $\WW(\bw)$ is quadratic such that $\partial^2 \WW(\bw) = H$, then the operator $A:\calX_c \rightarrow \calX^*$
  is linear and can be written in the triality form: $A  = \DD^* H \DD$,
which appears extensively in mathematical physics,     optimization, and linear systems \cite{gao-dual00,strang}.
Clearly, any convex quadratic function $\WW(\DD\uu)$ is objective due to the
 Cholesky decomposition $A = \Lam^* \Lam \succeq 0 $.

  In operations research,
    the  decision variable is  usually a vector   $\bu  \in \real^n$. If it
   represents  the products of a manufacture company,
   its dual variable $\barbu^* \in \real^n$ can be considered as market price (or demands), so the external energy
  $\FF(\bu) = - \la \bu , \barbu^*\ra =- \bu^T \barbu^* $ in this system
is  the total income of the company.
The products are produced by workers $\bw= \DD \bchi\in \real^m$ and
$\DD \in \real^{m\times n}$ is a matrix due to the cooperation.
 Workers are paid by salary $\bsig = \partial \WW(\bw)$,
therefore, the internal energy $\WW(\bw)$ in this example is the cost,
which could be an objective function (not necessary since the company is a man-made system).
Thus, $\Pi(\bu) = \WW(\DD \bu) + \FF(\bu)$ is the {\em total cost or target}  and the minimization problem $\min \Pi(\bchi)$ leads to
the equilibrium equation $
\DD^T \partial_{\beps} \WW(\DD \bu) = \barbu^*,
$
which is an algebraic equation in $\real^n$.

According to the action-reaction duality in physics,
if there is no action or demand (i.e. $\barbu^* = 0$), 
the system has non reaction (i.e. $\uu = 0$). 
Dually, for any given non-trivial input a real-world  problem  should have  at least one non-trivial solution.
\begin{definition}[Properly and Well-Posed Problems]
  A  problem   is called {\em properly posed} if
for any given non-trivial input it has  at least one non-trivial solution.
It is called {\em well-posed} if the solution is  unique. \vspace{-.2cm}
\end{definition}
Clearly, this definition is more general than Hadamard's well-posed problems in dynamical  systems
since   the continuity  condition is not required.
 Physically speaking, any real-world problems  should be well-posed since
  all natural phenomena  exist   uniquely.
  But practically, it is difficult to model  a real-world problem precisely.
  Therefore,  properly posed problems are allowed for the canonical duality theory.
  This definition is important for understanding the triality theory and NP-hard problems.\vspace{-.3cm}

\subsection{Nonconvex Analysis and Boundary-Value Problems}
For static systems, the unknown of a mixed boundary-value problem is  a vector-valued function
\[
\uu(\bx) \in \calU_a = \{ \bu\in \calC[\Oo, \real^m]| \;\; \bu(\bx) = \barbu \;\; \forall \bx \in \Gu\}, \;\;
\Oo\subset \real^d, \;\; d \le 3, \;\; m \ge 1, \;\; \partial \Oo = \Gu\cup\Gt,
\]
and the input is  $\barbu^* = \{ \bff(\bx) \;\;\forall  \bx \in \Oo, \;\; \bt(\bx) \; \forall \bx \in \Gt \}$ \cite{gao-strang89}.
In this case, the external energy is
$\FF(\bu) = -  \la \bu, \barbu^* \ra = - \int_\Oo \bu\cdot \bff \dO - \int_\Gt \bu\cdot \bt \dG$.
In nonlinear analysis, $\DD$ is a gradient-like partial differential operator and
$\ww = \DD \uu \in \calW_a \subset \calL^p[\Oo; \real^{m\times d}]$ is a
{\em two-point  tensor field} \cite{gao-dual00} over $\Oo$.
The internal energy $\WW(\ww)$ is defined by
\eb
\WW(\ww) = \int_\Oo \UU(\xx, \ww) \dO   ,
\ee
where $\UU(\xx,\ww):\Oo \times \calW_a \rightarrow \real$ is the {\em stored energy density}.
The system is  (space)  {\em homogeneous} if $\UU= \UU(\ww)$.
Thus, $\WW(\ww)$ is objective  if and only if $\UU(\xx,\ww)$ is  objective
on an objective set  $ \calW_a$.
By the facts  that $\ww = \DD \bu$ is a two-point tensor, which is not considered as a strain measure;
but   the (right) Cauchy-Green tensor $\bC = \ww^T \ww$ is an objective strain tensor,
there must exists a function  $\VV(\bC)$ such that $\WW( \ww) = \VV( \bC)$. In nonlinear elasticity,
the function $\VV(\bC)$ is usually convex and the duality $\bC^* = \partial\VV(\bC)$ is invertible (i.e. Hill's work-conjugate principle \cite{gao-dual00}).
These basic truths in continuum physics laid  a foundation for   the canonical duality theory.

By finite element method, the domain $\Oo$ is divided into $m$-elements $\{ \Oo^e \}$ such that
 the unknown function   is piecewisely  discretized by $\uu(\bx) \simeq \bN_e(\bx) \bchi_e \;\; \forall \bx \in \Oo^e$.
 Thus,   the nonconvex variational problem (\ref{eq-gs}) can be numerically reformulated in a
  global optimization problem
  \eb \label{eq-gop}
\min \{ \Pi(\bchi) = \WW(\bD \bchi) -  \la \bchi ,  \bff \ra  \;\; | \;\; \bchi \in \calX_c \} ,
\ee
where $\bchi = \{ \bchi_e \}$ is the discretized unknown $\bu(\bx)$, $\bD$ is a generalized (high-order) matrix  depending  on the interpolation $\bN_e(\bx)$ and  $\calX_c$ is
a  convex constraint set  including the boundary conditions.
The canonical dual finite element method was first proposed in 1996 \cite{gao-jem96}. Powerful applications
have been given recently in engineering and sciences \cite{gao-yu,santos-gao}.\vspace{-.3cm}

\subsection{Lagrangian Mechanics and Initial Value Problems}
   In Lagrange mechanics \cite{Lag,lan-lif}, the unknown $\uu(t) \in \calU_a \subset \calC^1[\Oo;\real^n]$
  is a vector field over a time domain   $\Oo \subset \real$.
  Its components $\{ u_i(t) \} \; (i = 1, \dots, n) $ are  known as
  the {\em Lagrangian coordinates}.  Its dual variable $\barbu^*$ is the action vector function in $\real^n$, denoted by  $\bff(t)$.
  The external energy $\FF(\uu) = - \la \uu, \barbu^* \ra = - \int_\Oo \uu(t) \cdot \bff(t) \dt$.
  While the internal energy $\WW(\DD\uu)$ is the so-called   action:
  \eb \label{eq-action}
  \WW(\DD \uu) = \int_\Oo L(t, \uu, \dot{ \uu} ) \dt , \;\;\;  L=  \TT( \dot{ \uu}  ) - \UU(t,  \uu),
  \ee
   where $\DD \uu = \{ 1 , \partial_t \} \uu = \{ \uu, \; \dot{\uu} \}$ is a vector-valued mapping,
  $\TT $ is the kinetic energy density,  $ \UU$ is the potential density, and  $L= \TT -  \UU$ is the
  {\em Lagrangian density}.    Together, $\Pi(\uu) = \WW(\DD \uu) + \FF(\uu)$   is called the {\em total action}.
  This standard  form  holds  from    the classical  Newton  mechanics to
  quantum field theory\footnote{See Wikipedia: \url{https://en.wikipedia.org/wiki/Lagrangian_mechanics}}.
  Its  stationary condition leads to the well-known {\em  Euler-Lagrange equation}:
   \eb
  A(\uu) =  \DD^* \partial\WW(\DD\uu) = \{ 1, - \partial_t \} \cdot \partial L(\uu, \dot{\uu}) =  - \partial_t  \partial_{\dot{\uu}}\TT( \dot{\uu})
     -  \partial_\uu \UU( t, \uu ) = \bff  .\label{eq-e-l}
   \ee
   The system is called (time) {\em homogeneous} if $L = L(\uu,\dot{\uu})$.
  In general, the kinetic energy $\TT $ must be   an objective function of the velocity\footnote{The objectivity of $\TT(\bv)$ is also called the isotropy in Lagrange  mechanics since $\bv $ is a vector (see \cite{lan-lif})}
   $\bv_k = \dot{\bx}_k(\uu)$
      of each particle $\bx_k = \bx_k(\uu) \in \real^3 \;\; \forall k \in \{1, \dots, K\}$,
   while the potential density $ \UU$ depends on each problem.  
For Newtonian mechanics, we have $\uu(t) = \bx(t)$ and $\TT(\bv) = \half  m \|\bv \|^2 $ is quadratic.
    If $ \UU = 0$, the equilibrium equation $A(\uu) = - m  \ddot{\bx}(t) = \bff$  includes
 the Newton  second law:  $\bF = m \ddot{\bx}$ and  the  third law:   $-\bF = \bff$.
 The first law  $\bv = \dot{\bx} = \bv_0$ holds only if $\bff = 0$.
  In this case, the system has either a trivial solution $\bx = 0$ or infinitely many solutions
  $\bx(t) = \bv_0 t + \bx_0$, depending on the initial conditions in $\calU_a$.
   This simple fact in elementary physics plays a key role in understanding the canonical duality theory and NP-hard problems in global optimization.

By using  the methods of  finite difference   and least squares \cite{gao-bridge,lat-gao-chaos},
the  general nonlinear dynamical system (\ref{eq-e-l}) can also be formulated as the  same global optimization problem
(\ref{eq-gop}),
where $\bchi = \{  u_i (t_k)\} $ is the Lagrangian coordinates $ u_i (i=1,\dots, n)$ at each discretized time  $t_k (k=1,\dots, m)$,   $\bD$ is a finite difference matrix   and $\calX_c$ is
a convex constraint set including the initial  condition \cite{lat-gao-chaos}.
By  the canonical duality theory, an intrinsic relation between  chaos in nonlinear dynamics and  NP-hardness in global optimization was revealed recently in   \cite{lat-gao-chaos}.



\subsection{Mono-/Bi-Dualities  and Duality Gap}
Lagrangian  duality was developed from  Lagrange mechanics since 1788 \cite{Lag}, where the kinetic energy $\TT(\bv) = \sum_{k}\half  m_k \|\bv_k\|^2$ is a  quadratic (objective) function.
For  convex static systems (or dynamical systems but $\UU(\bu) = 0$),
the stored energy $\WW:\calW_a \rightarrow \real$ is convex and its
Legendre conjugate   $\WW^*(\bsig)= \{\la \bw ; \bsig \ra - \WW(\bw ) | \;\; \bsig = \partial \WW(\bw) \} $ is uniquely defined
on $\calW^*_a$.
Thus,   by $\WW(\DD \uu) = \la \DD \uu; \bsig \ra - \WW^*(\bsig)$ the total potential $\Pi(\uu)$  can be   written
in the Lagrangian form\footnote{In Physics literature, the same notation $L$ is used  for both action
$L(\uu, \dot{\uu})$ and the Lagrangian ${L}(\uu, \bp)$ since both represent the same physical quantity.} $L:\calU_a \times \calW_a^* \rightarrow \real$:
\eb\label{eq-lagr}
L(\uu, \bsig) = \la \DD \uu ; \bsig \ra - W^*(\bsig) - \la \uu, \bff \ra
= \la \uu, \DD^* \bsig - \bff \ra - W^*(\bsig),
\ee
where $\uu\in \calU_a$ can be viewed as a Lagrange multiplier for the equilibrium equation $\DD^* \bsig = \bff \in \calU^*_a$.
In linear elasticity, $L(\uu, \bsig)$ is the well-known {\em Hellinger-Reissner complementary energy} \cite{gao-dual00}.
 Let
 $\calS_c = \{ \bsig \in \calW^*_a | \; \DD^* \bsig = \bff \} $ be the so-called
 {\em statically admissible space}. Then the Lagrangian dual  of the general  problem
   $(\calP_0)$ is  given by \cite{gao-dual00}
  \eb
(\calP_0^*): \;\;\;  \max \{  \Pi^*(\bsig) = - \WW^*(\bsig) | \; \bsig \in \calS_c \}, \label{eq-ld}
 \ee
and the saddle  Lagrangian leads to a well-known  min-max duality  in convex (static) systems:
\eb
\min_{\uu\in \calU_c} \Pi(\uu) = \min_{\uu\in \calU_a} \max_{\bsig \in \calW^*_a} L(\uu, \bsig) =
  \max_{\bsig \in \calW^*_a} \min_{\uu\in \calU_a} L(\uu, \bsig) =
 \max_{\bsig \in \calS_c} \Pi^*(\bsig)  .
\ee
This one-to-one duality is the so-called the  {\em mono-duality} in Chapter 1 \cite{gao-dual00}, or
 the  {\em complementary-dual variational principle} in continuum physics. 
In finite  elasticity,  the Lagrangian dual is also known as the {\em Levison-Zubov principle}.
However, this principle holds only for convex problems.

For convex Hamiltonian systems the action  $\WW(\DD\bu) $ in (\ref{eq-action}) is a d.c.  (difference of convex) functional   and the Lagrangian
has its standard form in Lagrangian mechanics (see Chapter 2.5.2   \cite{gao-dual00}
with $\bu = \bq(t)$ and $\bsig = \bp$):\vspace{-.3cm}
 \eb \label{eq-lag}
 L(\bq, \bp) = \la\dot{\bq} ; \bp \ra - \int_\Oo [\TT^*(\bp) + \UU(\bq)]\dt - \la \bq , \bff \ra,\vspace{-.2cm}
 \ee
where $\bq \in \calU_a \subset \calC^1[\Oo,\real^n]$ is the Lagrange coordinate and
$\bp \in \calS_a  \subset \calC[\Oo,\real^n]$ is the momentum. In this case,  the Lagrangian is a bi-concave functional on $\calU_a \times \calS_a$, but the Hamiltonian $H(\bq,\bp) = \la\DD{\bq} ; \bp \ra  -\LL(\bq,\bp) $ is convex\footnote{This is the reason that  instead of the Lagrangian, the Hamiltonian is extensively used in dynamics.}.
The total action and its canonical dual are \cite{gao-dual00} \vspace{-.2cm}
\begin{eqnarray}
\Pi(\bq)& = & \max \{ \LL(\bq, \bp) | \;\; \bp \in \calV^*_a\} = \int_\Oo[\TT(\dot{\bq}) - \UU(\bq) ] \dt - \la \bq, \bff \ra \;\; \forall \bq \in \calU_c \\
\Pi^d(\bp) &  = & \max \{ \LL(\bq, \bp) | \;\; \bq \in \calU_a\} = \int_\Oo[\UU^*(\dot{\bp}) - \TT^*(\bp) ] \dt  \;\; \forall \bp \in \calS_c\vspace{-.3cm}
\end{eqnarray}
Although  both of them are d.c.  functionals,
the duality between the kinetic energy $\TT(\dot{\bq})$ and the potential $\UU(\bq)$
 leads to a so-called {\em bi-duality  }  first presented in author's book
 Chapter 2 \cite{gao-dual00}:
 \eb
 \min \Pi(\bq) = \min \Pi^d(\bp), \;\; \max \Pi(\bq) = \max \Pi^d(\bp).
 \ee
 The mathematical proofs of this  theory   were  given in Chapter 2.6 \cite{gao-dual00}  for  convex Hamiltonian systems  and  in Corollary 5.3.6 \cite{gao-dual00} for nonconvex  programming problems.
 This bi-duality revealed not only an interesting dynamical extremum principle in periodic motion, but also
 an important truth in convex Hamiltonian systems (see page 77 \cite{gao-dual00}):
 {\em the least action principle is incorrect for any periodic motion, it holds only for linear potential $\UU(\bq)$}\footnote{
This truth is not known to many people in physics. Only in a  footnote  of the celebrated book  (Section 1.2, \cite{lan-lif})
Landau and  Lifshitz   pointed out
 that the least action principle  holds only for a sufficient small of time interval, not for
 the whole trajectory of the system.}.

In real-world problems the stored energy $\WW(\bw)$ is usually  nonconvex in order to model complex phenomena. Its complementary energy can't be determined uniquely by the Legendre transformation.
Although its  Fenchel conjugate $\WW^\sharp:\calW^*_a \rightarrow \real \cup \{ + \infty\}$ can be uniquely defined,
the Fenchel-Moreau  dual problem
\eb
(\calP^\sharp_0): \;\;\;
\max \{ \Pi^\sharp(\bsig) = - \WW^\sharp(\bsig) | \;\; \bsig \in \calS_c\}
\ee
is not considered as a complementary-dual problem due to Fenchel-Young inequality:
\eb
\min \{  \Pi(\uu)  | \; \uu\in \calU_c \} \ge \max  \{ \Pi^\sharp(\bsig)  | \;\; \bsig \in \calS_c\} ,
\ee
and $\theta = \min  \Pi(\uu) - \max \Pi^\sharp (\bsig) \neq 0$ is the well-known  {\em duality gap}.
This duality gap is intrinsic to all  Lagrange-Fenchel-Moreau types duality problems since  the
linear operator $\DD$ can't change the nonconvexity of $\WW(\DD\uu)$.
It turns out that the existence of a pure stress based complementary-dual principle  has been a well-known debate in nonlinear  elasticity
for more than fifty years \cite{li-gupta}.

\begin{remark}[Lagrange Multiplier Law]
{\em
Strictly  speaking, the Lagrange multiplier method can be used mainly for equilibrium constraint
 in $\calS_c$  and  the Lagrange multiplier must be the solution to the primal problem (see Section 1.5.2 \cite{gao-dual00}).
 The equilibrium equation $\DD^* \bsig = \bff$
  must be an invariant under certain coordinates transformation,
say the law of angular momentum conservation,
 which is  guaranteed by  the objectivity of  the stored energy $\WW(\DD\uu)$ in continuum mechanics (see Definition 6.1.2, \cite{gao-dual00}), or by the
 isotropy  of the kinetic  energy  $\TT(\dot{\uu})$ in Lagrangian mechanics \cite{lan-lif}.
 Specifically,  the equilibrium equation for Newtonian  mechanics is an invariant under the
 {\em  Calilean transformation};
 while for Einstein's special relativity theory, the  equilibrium equation $\DD^* \bsig = \bff$ is an invariant under the
 {\em Lorentz transformation.}
 For linear equilibrium equation, the  quadratic
 $\WW(\bw)$ is naturally an objective function for convex systems.
    Unfortunately, since  the concept of the  objectivity is misused in mathematical optimization and the notation of the Euclidian coordinate $\bx = \{ x_i\}$ is  used as the unknown,
    the Lagrange multiplier method and augmented methods have been mistakenly used for solving general nonconvex optimization problems, which produces
    many  artificial  duality gaps  \cite{lat-gao-opl}.
     }
     \end{remark}
\section{Unified Problem and Canonical Duality-Triality Theory}
In this section, we simply restrict our discussion  in
  finite-dimensional space $\calX$. Its element $\bchi \in \calX$ could be a vector,  a  matrix, or a tensor\footnote{Tensor is a geometrical object in mathematics and physics,   which is defined as a multi-dimensional array satisfying a transformation law, see \url{https://en.wikipedia.org/wiki/Tensor}. A tensor must be independent of a particular choice of coordinate system (frame-invariance).
  But this  terminology has been also misused  recent years in optimization literature such that any multi-dimensional array of data
  is called tensor.}. In this case, the linear operator $\DD$ is a generalized matrix\footnote{A generalized matrix
  $\bD = \{ D^{i \cdots j}_{\alp \cdots \gamma}\}$ is a multi-dimensional array but not necessary to satisfy  a transformation law, so it is not a tensor. In order to avoid confusion, it can be called a {\em tentrix}.}
    $\bD:\calX  \rightarrow \calW$ and $\calW$ is a generalized matrix space equipped  with
   a  natural norm $\|\bw\|$.
  Let $\calX_a \subset \calX$ be a convex subset  and $\calX^*_a$ be its dual set such that
  for any given input $\bff \in \calX^*_a$   the subjective function $   \la \bchi, \bff \ra \ge 0 \;\; \forall \bchi \in \calX_a$.  Then  the
    multi-scale optimization problem (\ref{eq-gop}) can be re-proposed  as
   \eb\label{eq-unimod}
   (\calP): \;\;
  \min \left\{ \Pi(\bchi) = \WW(\bD\bchi)  - \la \bchi , \bff \ra | \;\; \bchi \in \calX_c  \right\},
  \ee
where $\calX_c  = \{ \bchi \in \calX_a | \;\; \bD \bchi  \in \calW_a \}$.
  Although the objectivity is necessary for  real-world modeling, the numerical discretization of $\WW(\DD\uu)$
  could lead to a complicated function $\WW(\bD \bchi)$, which may not be  objective in $\bw = \bD \bchi$.
Also in operations research, many challenging problems are artificially proposed.
Thus, the objectivity required in  Gao and Strang's work  on nonlinear elasticity has been  relaxed by the canonical  duality  since  2000 \cite{gao-jogo00}.


\subsection{Canonical Transformation and Gap Function}


In the canonical duality theory,
a real-valued function
$\VV:\calE_a\subset \calE  \rightarrow \real$ is said to be {\em canonical} if the duality relation
$ \bvsig  = \partial \VV(\bxi):\calE_a \rightarrow \calE^*_a \subset \calE^*$ is bijective.
The canonical duality is a fundamental principle
 in sciences  and oriental philosophy, which underlies  all natural phenomena. Therefore, instead of the objectivity in continuum physics,
a generalized  objective function $\WW(\bw)$ is used in the canonical duality theory under  the following assumption.
 \begin{assumption}
For a given  $\WW:\calW_a \rightarrow \real$,
there exists a canonical  measure $\bxi:\calW_a \rightarrow \calE_a$ and  a canonical function $\VV:\calE_a \rightarrow \real$ such that the following conditions hold:

(A1.1) Positivity: $\WW(\bw) \ge 0 \;\; \forall \bw \in  \calW_a$;

(A1.2) Canonicality:     $\WW(\bw) = \VV(\bxi(\bw)) \;\; \forall \bw \in \calW_a$; And either

(A1.3) Coercivity: $\lim \WW(\ww) = \infty $ as $\| \ww \| \rightarrow \infty$,  or

(A1.3*) Boundness: $\calW_a$ is bounded.
\end{assumption}
Generally speaking, the conditions (A1.1) and (A1.2) are necessary for any real-world problems;
while (A1.3) and (A1.3*) depend  mainly on the magnitude of the input $\bff \in \calX^*_a$.
Usually,  the  coercivity  is for  small $\|\bff\|$ (within the system's capacity, such as elasticity)  and  the boundness  is for big  $\|\bff\|$ (beyond the
system's capacity, such as plasticity) \cite{gao-jmaa98}.

Let $ \Lam = \bxi  \circ \bD :\calX_a \rightarrow \calE_a$ be the so-called  {\em geometrically admissible operator}.
 The canonicality $\WW(\bD \bchi) = \VV(\Lam(\bchi))$ is also called the {\em canonical transformation} in the canonical duality theory.
Let $\la \bxi ;  \bvsig \ra :\calE \times \calE^* \rightarrow \real$ be the bilinear form  which puts $\calE$ and $\calE^*$ in duality. By  (A1.2),  we have  $\calX_c  = \{ \bchi \in \calX_a | \;\; \Lam(\bchi) \in \calE_a \}$
  and the problem $(\calP)$
    can be equivalently reformulated in the following canonical form \vspace{-.2cm}
   \eb
   (\calP): \;\;
  \min \left\{ \Pi(\bchi) = \VV(\Lam( \bchi)) - \la \bchi , \bff \ra | \;\; \bchi \in \calX_c  \right\}.
  \ee
  By the facts that the canonical  duality is a universal principle  in nature  and the canonical measure
  $\Lam(\bchi)$ is not necessarily to be objective,
  the canonical transformation holds
for   general  problems and  the problem $(\calP)$ can be used to model general  complex systems.
The criticality condition of $(\calP)$ is governed by the {\em fundamental principle of virtual work}:\vspace{-.2cm}
\eb
\la \Lam_t (\bchi)\delta \bchi ;  \bvsig \ra = \la \delta  \bchi , \Lam^*_t(\bchi)  \bvsig  \ra = \la \delta \bchi, \bff \ra  \;\; \forall \delta \bchi\in \calX_c,
\ee
where $\Lam_t(\bchi) = \partial \Lam(\bchi)$ represents a generalized  \G (or directional) derivative of $\Lam(\bchi)$,
its adjoint   $\Lam^*_t$ is called the  balance operator,   $ \bvsig= \bC^*(\bxi)  = \partial \VV(\bxi)$
  and $\bC^*:\calE_a \rightarrow \calE^*_a $ is a canonical dual (or constitutive) operator.
The strong form of this   virtual work  principle  is called {\em
the canonical equilibrium equation}:
\eb
\bA (\bchi) = \Lam^*_t(\bchi) \bC^* (\Lam(\bchi)) = \bff .
\ee
A system governed by this equation is  called
a {\em canonical system} and is denoted as (see Chapter 4, \cite{gao-dual00})
\[
\mathbb{S}_a = \{ \la \calX  , \calX^* \ra ,   \la  \calE  ;   \calE^* \ra  ;  ( \Lam,  \bC^* )  \}.
\]
\begin{definition}[Classification of Nonlinearities]
 The system $\mathbb{S}_a $ is called {\em geometrically nonlinear} (resp. linear) if the geometrical operator
$\Lam:\calX_a \rightarrow \calE_a$ is nonlinear (resp. linear);
the system  is called {\em physically (or constitutively) nonlinear } (resp. linear)
if the canonical dual operator $\bC^*:\calE_a \rightarrow \calE^*_a $ is nonlinear (resp. linear);
the systems is called {\em fully nonlinear} (resp. linear) if
it is  both geometrically and physically nonlinear (resp. linear).
\end{definition}
Both geometrical and physics nonlinearities are basic concepts in  nonlinear field theory.
The mathematical definition was first given by the author in 2000 under the canonical transformation \cite{gao-jogo00}.
  A diagrammatic representation of this canonical system  is
shown in  Figure \ref{frame1}.

\begin{figure}[h]
\begin{center}
\setlength{\unitlength}{.5mm}
\begin{picture}(170,65)(50,5)
 \put(125,56){\makebox(0,0){$ \bA $}}
  \put(125,44){\makebox(0,0){$ \la \bchi  \bif \bchi^* \ra $}}
 \put(97,50){\vector(1,0){55}}

 \put(97,10){\vector(1,0){55}}
  \put(125,3){\makebox(0,0){$ \bC^*$}}
 \put(125,15){\makebox(0,0){$ \la \bxi  \bis    \bvsig  \ra $}}

\put(86,30){\makebox(0,0)[r]{$  \Lam_t + \Lam_c =  \Lam $  }}
\put(90,42){\vector(0,-1){24}}
\put(166,30){\makebox(0,0)[l] {$\;\; \Lam^*_t = (\Lam - \Lam_c)^* $}}
\put(160,18){\vector(0,1){24}}
 \put(90,10){\makebox(0,0){$\calE_a$}}
\put(77,10){\makebox(0,0){$   \bxi \in  $}}

\put(90,50){\makebox(0,0){$\calX_a$}}
\put(77,50){\makebox(0,0){$ \bchi\in   $}}
 \put(160,10){\makebox(0,0){$\;\;\calE^*_a $}}
\put(176,10){\makebox(0,0){$  \ni   \bvsig $}}
\put(160,50){\makebox(0,0){ $\;\;\calX^*_a$}}
\put(178,50){\makebox(0,0){$  \ni \bchi^*  $}}
\end{picture}
\caption{Diagrammatic representation  for a canonical   system }\label{frame1}
\end{center}
\end{figure}
This diagram shows a symmetry broken in the canonical equilibrium equation,
i.e., instead of $\Lam^*$, the balance operator  $\Lam_t^*$ is adjoined with $\Lam_t$.  It was discovered by Gao and Strang \cite{gao-strang89}
 that by introducing a complementary operator $\Lam_c(\bchi) = \Lam(\bchi) - \Lam_t(\bchi) \bchi$,  this locally broken symmetry is  recovered by a so-called  complementary gap function
 \eb
 \Gap(\bchi, \bvsig) = \la - \Lam_c(\bchi) ; \bvsig \ra,
 \ee
  which plays a key role
in global optimization  and the triality theory.
Clearly, if $\Lam = \bD$ is linear, then $\Gap = 0$. Thus, the following statement is important to understand complexity:

 {\em Only  the geometrical  nonlinearity  leads to nonconvexity in optimization, bifurcation in analysis, chaos in dynamics, and NP-hard problems in complex systems}.

\subsection{Complementary-Dual Principle and Analytical Solution}
For a given canonical function
$\VV:\calE_a \rightarrow \real$, its conjugate
 $\VV^*:\calE^*_a \rightarrow \real$ can be uniquely defined by the Legendre transforation \vspace{-.3cm}
\eb
\VV^*( \bvsig ) = \sta \{ \la \bxi ;  \bvsig  \ra - \VV(\bxi) |  \;\;  \bxi  \in \calE_a \},
\ee
where $\sta \{ f(\bchi) | \; \bchi \in \calX \}$ stands for finding the stationary value of $f(\bchi)$
on $\calX$, and  the following  canonical  duality relations  hold on $\calE_a \times \calE^*_a$:
\eb \label{eq-cdr}
 \bvsig =  \partial     \VV(\bveps) \;\;\Leftrightarrow \;\; \bveps = \partial   %
\VV^*(\bvsig) \;\;\Leftrightarrow \;\; \VV(\bveps) + \VV^*(\bvsig) =  \la  \bveps ;  \bvsig \ra.
\ee

If the canonical function is  convex and lower semi-continuous, the \G derivative $\partial$ should be replaced by  the sub-differential and  $\Phi^*$ is replaced by the Fenchel conjugate $\Phi^\sharp(\bvsig) = \sup\{ \la \bxi ; \bvsig\ra - \Phi(\bxi) | \; \bxi \in \calE_a \}$. In this case,  (\ref{eq-cdr}) is replaced by the generalized canonical duality
\eb \label{eq-gcdr}
 \bvsig \in  \partial     \VV(\bveps) \;\;\Leftrightarrow \;\; \bveps \in \partial   %
\VV^\sharp(\bvsig) \;\;\Leftrightarrow \;\; \VV(\bveps) + \VV^\sharp(\bvsig) =  \la  \bveps ;  \bvsig \ra\;\; \forall (\bxi, \bvsig) \in \calE_a \times \calE^*_a .
\ee
If the  convex set $\calE_a$ contains   inequality constrains, then (\ref{eq-gcdr}) includes
all the {\em internal  KKT conditions} \cite{lat-gao-opl,gao-jmaa98}. In this sense, a
 KKT point  of the canonical form $\Pi(\bchi)$ is a generalized  critical point of $\Pi(\bchi)$.

By the complementarity  $\VV(\Lam(\bchi)) =    \la  \Lam(\bchi) ;  \bvsig \ra  - \VV^*(\bvsig) $,
the canonical form of $\Pi(\bchi)$ 
can be equivalently written in
  Gao and Strang's   {\em total complementary function} $\Xi:\calX_a \times \calE_a^* \rightarrow \real $   \cite{gao-strang89}
\eb
\Xi(\bchi, \bvsig) =   \la \Lam(\bchi) ;  \bvsig \ra   - \VV^*(\bvsig)
  - \la \bchi, \bff \ra   .
\ee
Then, the canonical dual function $\Pi^d:\calS_c \rightarrow \real$ can be obtained by
the  {\em canonical dual transformation}:
\eb
\Pi^d(\bchi) = \sta \{ \Xi(\bchi, \bvsig) | \;\; \bchi \in \calX_a \} = \Gap^\Lam(\bvsig) - \Phi^*(\bvsig),
\ee
where $\Gap^\Lam (\bvsig) = \sta \{ \la \Lam(\bchi ) ;  \bvsig \ra - \la \bchi, \bff \ra | \;\;  \bchi \in \calX_a \}$,
which is  defined on the canonical dual feasible space $\calS_c = \{ \bvsig \in \calE^*_a | \; \;  \Lam_t^*(\bchi)  \bvsig =  \bff  \;\; \forall \bchi \in \calX_a \}.$ Clearly, $\calS_c \neq \emptyset$ if $(\calP)$ is properly posed.
\begin{thm}[Complementary-Dual Principle \cite{gao-dual00}]\label{thm-comd}
The pair  $(\barbchi, \barbvsig)$ is a critical point
   of $\Xi(\bchi, \bvsig)$ if and only if $\barbchi$ is a critical point  of $\Pi(\bchi)$ and
  $\barbvsig$ is a critical point of $\Pi^d(\bvsig)$. Moreover,
\eb
\Pi(\barbchi) = \Xi(\barbchi, \barbvsig) = \Pi^d(\barbvsig).
\ee
\end{thm}
{\bf Proof}. The criticality condition $\partial\Xi(\barbchi, \barbvsig) = 0$ leads to the following canonical equations
\eb
 \Lam(\barbchi) = \partial \VV^*(\barbvsig), \;\; \;  \Lam_t^*(\barbchi) \barbvsig  = \bff.
\ee
The theorem is proved by  the canonical duality (\ref{eq-cdr}) and the definition of $\Pi^d$. \hfill $\Box$

Theorem \ref{thm-comd} shows a one-to-one correspondence of the critical points between the primal function and its canonical dual. In large deformation theory, this theorem solved the  fifty-year   old open problem on
  complementary variational principle
and is known as the Gao principle in literature \cite{li-gupta}.

In  real-world applications, the geometrical operator $\Lam$ is usually  {\em quadratic homogeneous}
i.e. $\Lam (\alp \bchi) = \alp^2 \Lam(\bchi) \;\; \forall \alp \in \real$. In this case, we have  \cite{gao-strang89}
\[
\Lam_t (\bchi) \bchi = 2 \Lam(\bchi), \;\; \Lam_c(\bchi) = - \Lam(\bchi) , \mbox{ and}
\]
\eb
\Xi(\bchi, \bvsig) =   \Gap(\bchi, \bvsig)   - \VV^*(\bvsig)
  - \la \bchi, \bff \ra = \half \la \bchi, \bG(\bvsig) \bchi \ra - \VV^*(\bvsig) - \la \bchi,  \bff  \ra  ,  \label{eq-tcf}
\ee
where $ \bG(\bvsig) = \partial^2_{\bchi} \Gap(\bchi, \bvsig)$.
 Then,    the canonical dual  function  $\Pi^d(\bvsig)$ can written explicitly as
\eb
\Pi^d(\bvsig)  = \{ \Xi(\bchi, \bvsig) | \;\; \bG(\bvsig) \bchi = \bff  \;   \forall \bchi \in \calX_a\}
= - \half \la [\bG(\bvsig)]^{+} \bff  , \bff \ra - \VV^*(\bvsig),
\ee
where $\bG^+$ represents a generalized inverse of $\bG$.

\begin{thm}[Analytical Solution Form\cite{gao-dual00}]
If $\barbvsig \in \calS_c$ is a critical point of $\Pi^d(\bvsig)$, then
\eb \label{eq-ana}
\barbchi = [\bG(\barbvsig)]^{+} \bff
\ee
is a critical point  of $\Pi(\bchi)$ and $
\Pi(\barbchi) = \Xi(\barbchi, \barbvsig) = \Pi^d(\barbvsig).
$ Dually, if $\barbchi\in \calX_c$ is a critical point of $\Pi(\bchi)$, it must be in the form of (\ref{eq-ana})
for a critical point $\barbvsig\in \calS_c$ of $\Pi^d(\bvsig)$.
\end{thm}

This unified analytical solution form holds not only for general global optimization problems in finite dimensional systems,
but also for a large-class of  nonlinear boundary/initial value problems
 in nonconvex analysis   and dynamic systems \cite{gao-advances}.

\subsection{Triality Theory and NP-Hard Criterion}
\begin{definition}
[Degenerate and Non-Degenerate Critical Points,  Morse Function]$\;$\hfill
Let $\barbchi \in \calX_c$ be a critical point of a real-valued function $\Pi:\calX_c \rightarrow \real$. $\barbchi $ is  called  degenerate (res. non-degenerate)  if
 the Hessian matrix of $\Pi(\bchi)$  is singular (resp. non-singular) at $\barbchi $.
The function $\Pi:\calX_c \rightarrow \real$ is called a Morse function if it has no degenerate critical points.
\end{definition}

\begin{thm}[Triality Theory \cite{gao-jogo00}]\label{thm-tri}
Suppose that  
 $\VV:\calE_a \rightarrow \real$ is convex,   $(\barbchi, \barbvsig)$ is a non-degenerate critical point
   of $\Xi(\bchi, \bvsig)$  and   $\calX_{o}\times \calS_{o} $  is  a neighborhood \footnote{The neighborhood $\calX_o$ of $\barbchi$  means that on which, $\barbchi$ is the only critical  point (see page 140 \cite{gao-dual00}).} of $\left(
\barbchi,\barbvsig \right) $.

 If  $\barbvsig \in \calS^+_c = \{ \bvsig \in \calS_c | \;\; \bG(\bvsig )  \succeq 0 \}$,
 then
\begin{equation}
\Pi(\barbchi) = \min_{\bchi \in \calX_c}  \Pi\left( \bchi\right)  =\max_{\bvsig \in \calS^+_c}  \Pi^{d}\left(\bvsig \right) = \Pi^d(\barbvsig) . \label{Global0}
\end{equation}

If $\barbvsig \in \calS^-_c = \{ \bvsig \in \calS_c| \;\; \bG(\bvsig )  \prec 0 \}$,
then  we have either    
\begin{equation}
\Pi(\barbchi) = \max_{\bchi \in \calX_o}  \Pi\left( \bchi\right)  =\max_{\bvsig \in \calS_o}  \Pi^{d}\left(\bvsig \right) = \Pi^d(\barbvsig)  ,  \label{dmax}
\end{equation}
or (if $\dim \Pi = \dim \Pi^d$)    
\begin{equation}
\Pi(\barbchi) = \min_{\bchi \in \calX_o}  \Pi\left( \bchi\right)  =\min_{\bvsig \in \calS_o}  \Pi^{d}\left(\bvsig \right) = \Pi^d(\barbvsig) .  \label{dmin}
\end{equation}
\end{thm}

The statement (\ref{Global0}) is the so-called {\em canonical min-max duality}, which can be proved easily by  Gao and Strang's work in 1989 \cite{gao-strang89}.
Clearly,  $\bvsig \in \calS^+_c$ if and only if $\Gap(\bchi, \bvsig) \ge 0 \;\; \forall \bchi \in \calX$.
This duality theory shows that the Gao-Strang gap function
  provides a    global optimum criterion.
   The statements (\ref{dmax}) and (\ref{dmin}) are called the  {\em canonical double-max} and   {\em double-min dualities}, respectively,  which can be used to find local extremum solutions.

  The triality theory shows that the nonconvex minimization problem $(\calP)$ is canonically dual
   to  the following maximum  stationary  problem
   \eb \label{eq-pidcor}
 (\calP^d): \;\; \max \sta \{ \Pi^d(\bvsig) | \;\; \bvsig \in \calS^+_c \}.
 \ee
 \begin{thm}[Existence and Uniqueness Criteria \cite{gao-cace}]
 For a properly posed  $(\calP)$, if the canonical function $\Phi:\calE_a \rightarrow \real$ is convex, {\em int}$ \calS^+_c \neq \emptyset$,  and 
\eb
\lim_{\alp \rightarrow 0^+}  \Pi^d(\bvsig_0 + \alp \bvsig) = - \infty \;\; \forall \bvsig_o \in \partial \calS^+_c, \;\; \forall \bvsig \in \calS^+_c,
 \ee
 then $(\calP^d)$  has at least one solution $\barbvsig\in\calS^+_c$
 and $\barbchi = [\bG(\barbvsig) ]^+\bff$ is a solution to $(\calP)$.
The solution is unique if $\bH = \partial \bG(\barbvsig) \succ 0$.
\end{thm}
{\bf Proof}.
Under the required conditions 
   $-\Pi^d:\calS^+_c \rightarrow \real$ is convex and coercive  and  int$\calS^+_c \neq \emptyset$. 
 Therefore, $(\calP^d)$ has at least one  solution. If $\bH \succ 0$,  then  $\Pi^d:\calS^+_c \rightarrow \real$
  is strictly concave and $(\calP^d)$ has a unique solution.\hfill $\Box$

 This theorem shows that  if    int$\calS^+_c \neq \emptyset$ 
 the  nonconvex problem  $(\calP)$ is canonically  dual to $(\calP^d)$ which can be solved easily.
 Otherwise, the problem $(\calP)$ is canonically dual to the following minimal stationary problem, i.e. to find a global
 minimum stationary value of $\Pi^d$ on $\calS_c$:
 \eb
 (\calP^s): \;\; \min \sta \{ \Pi^d(\bvsig) | \;\; \bvsig \in \calS_c \},
 \ee
   which could be really NP-hard since  $\Pi^d(\bvsig)$ is nonconvex on the nonconvex set $\calS_c$.
   Therefore, a conjecture was proposed in \cite{gao-jimo07}.
   \begin{Conj}[Criterion of NP-Hardness]
  A properly posed problem $(\calP)$ is NP-hard if and only if   {\em int }$\calS^+_c = \emptyset$.
  \end{Conj}

The triality theory was  discovered by the author
 during his research on post-buckling of a large deformed elastic beam in 1996 \cite{gao-amr97}, where the primal variable $\bu(\bx)  $ is a displacement vector in $\real^2$ and   $\bvsig(\bx)$ is a canonical dual stress also in $\real^2$.
Therefore, the triality theory was correctly proposed in nonconvex analysis, which provides for the first time
a complete set of solutions to  the  post-buckling problem. 
Physically, the  global  minimizer $\barbu(\bx)$ represents a stable
 buckled beam configuration (happened naturally),
the local minimizer is an unstable buckled state (happened occationally),
while the local maximizer is  the unbuckled beam state.
Mathematical proof of the triality  theory was given in \cite{gao-dual00}
for one-D nonconvex variational problems  (Theorem 2.6.2)
and for finite dimensional optimization problems (Theorem 5.3.6 and Corollary 5.3.1).
 In 2002, the author  discovered some countexamples to the canonical double-min duality when $\dim \Pi \neq \dim \Pi^d$ and this statement was removed from the triality theory (see Remark 1 in \cite{gao-opt03} and Remark for Theorem 3 in \cite{gao-amma03}).
 Recently, the author and his co-workers proved that the canonical double-min duality holds weakly
 when $\dim \Pi \neq \dim \Pi^d$ \cite{gao-wu,chen-gao-jogo,mora-gao-memo}.
 It was also discovered by using the canonical dual finite element method that the local minimum solutions  in nonconvex mechanics  are very sensitive not only to the input and boundary conditions of a given system, but also to such artificial conditions as
 the  numerical discretization and  computational precision, etc.
 The triality theory provides a precise mathematical tool for studying and understanding complicated natural phenomena.

The canonical duality-triality theory   has
been successfully used for solving a wide class problems in both global optimization
and  nonconvex  analysis \cite{gao-bridge}, including certain challenging
 problems  in  nonlinear PDEs and large deformation mechanics  \cite{gao-anti}.\vspace{-.3cm}

\section{ Applications in Complex Systems}
Applications to nonconvex constrained global optimization have been discussed in \cite{lat-gao-opl,g-r-s}.
This section presents applications to two general    global optimization problems.\vspace{-.3cm}

\subsection{Unconstrained Nonconvex Optimization Problem}  \vspace{-.3cm} 
\eb
(\calP): \;\;\; \min \left\{ \Pi(\bchi)  = \sum_{s=1}^m \Phi_s(\Lam_s(\bchi)) - \la  \bchi, \bff \ra | \;\; \bchi\in \calX_c\right\} ,
\ee
where the canonical measures  $\bxi_s= \Lam_s(\bchi)$ could be either a scalar 
or a generalized matrix,
$\Phi_k(\bxi_k) $ are any  given canonical functions, such as
polynomial, exponential, logarithm, and their compositions, etc.
For example, if $\bchi \in  \calX_c \subset \real^n$ and  \vspace{-.3cm}
\begin{eqnarray}
\WW(\bD \bchi) &=& \sum_{i\in \mathbb{I}}  \half {\alp_i} \bchi^T \bQ_i \bchi
 + \sum_{j \in \mathbb{J }}   \half \alp_j \left( \half \bchi^T \bQ_j \bchi +  \beta_j \right)^2 \nonumber \\
& +&  \sum_{k \in \mathbb{K }}  \alp_k \exp \left( \half   \bchi^T \bQ_k \bchi    \right)
+ \sum_{\ell  \in \mathbb{L}} \half \alp_\ell   \bchi^T \bQ_\ell\bchi   \log (\half  \bchi^T \bQ_\ell\bchi )  , \vspace{-.3cm}
\end{eqnarray}
where $\{ \bQ_s \} $ are positive-definite matrices to allow the Cholesky decomposition
 $\bQ_s = \bD^T_s \bD_s$ for all
$ s \in \{ \mathbb{I}, \mathbb{J}, \mathbb{K}, \mathbb{L}\}$ and $\{ \alp_s, \beta_s \} $ are physical constants, which could be either positive or negative
under Assumption 1.
This  general function includes naturally the so-called d.c. functions (i.e. difference of convex functions).
By using the canonical measure
\[
\bxi = \{ \xi_s \}  = \left \{\half {\alp_i} \bchi^T \bQ_i \bchi ,  \half \bchi^T \bQ_r \bchi \right \} \in \calE_a = \real^p \times \real^{q}_+ , \;\; p = \dim \mathbb{I}, \;\; q = \dim\mathbb{J} + \dim\mathbb{K}+ \dim\mathbb{L}
\]
 where    $\real^q_+ = \{ \bx \in \real^q| \; x_i \ge 0 \; \forall i = 1, \dots, q\}$,
 $\WW(\bw)$  can be written in the canonical form
 \begin{eqnarray*}
    \VV (\bxi)  &=&  \sum_{i\in \mathbb{I}}  \xi_i +  \sum_{j\in \mathbb{J}} \half \alp_j (\xi_j + \beta_j)^2 +
 \sum_{k\in \mathbb{K}} \alp_k \exp \xi_k   + \sum_{\ell\in \mathbb{L}} \alp_\ell  \xi_\ell \log\xi_\ell  .
 \end{eqnarray*}
Thus, $  \partial \Phi(\bxi) = \{ 1,  \vsig_r\} $ in which,
$\bvsig = \{ \alp_j (\xi_j + \beta_j), \alp_k \exp \xi_k, \; \alp_\ell ( \log\xi_\ell - 1) \} \in \calE^*_a$ and
\[
\calE^*_a  =
\{ \bvsig\in \real^q |\; \vsig_j \ge - \alp_j \beta_j \;\; \forall j \in \mathbb{J}, \;\; \vsig_k \ge \alp_k
\;\;\forall k \in \mathbb{K}, \;\;
\vsig_\ell \in \real \;\; \forall \ell \in \mathbb{L}\}.
\]
The conjugate of $\Phi$   can be easily obtained as
\begin{eqnarray}
\Phi^*(\bvsig)  &= &   \sum_{j\in \mathbb{J}}  \left(\frac{1}{2 \alp_j} \vsig_j^2 + \beta_j \vsig_j \right)
+ \sum_{k\in \mathbb{K}}   \vsig_k (\ln(\alp_k^{-1} \vsig_k) - 1)
+ \sum_{\ell\in \mathbb{L}}\alp_\ell \exp(\alp_\ell^{-1} \vsig_\ell -1).
\end{eqnarray}
Since $\Lam(\bchi)$ is  quadratic homogenous, the gap function $\Gap$  and its $\Lam$-conjugate $\Gap^\Lam$ in this case are
\[
\Gap(\bchi, \bvsig) = \half \bchi^T \bG(\bvsig) \bchi, \;\;
\Gap^\Lam(\bvsig) = \half \bff^T [\bG(\bvsig)]^+ \bff, \;\;
\bG(\bvsig)  = \sum_{i \in \mathbb{I}}   \alp_i  \bQ_i
+ \sum_{s\in \{ \mathbb{J}, \mathbb{K}, \mathbb{L}\}}   \vsig_s  \bQ_s    .
\]
Since $\Pi^d(\bvsig) = -\Gap^\Lam(\bvsig) - \Phi^*(\bvsig) $ is concave and  $\calS^+_c$ is a  closed convex set,
if for the given physical constants and the input $\bff$ such that $ \calS^+_c\neq \emptyset$, the canonical dual  problem $(\calP^d)$
has at least one solution $\barbvsig\in \calS^+_c \subset \real^q$ and
$\barbchi = [\bG(\barbvsig)]^+ \bff \in\calX_c \subset \real^n$ is a global minimum solution to $(\calP)$.
If  $n \gg q$, the problem $(\calP^d)$ can be much easier than $(\calP)$.

\subsection{Mixed Integer Nonlinear Programming (MINLP)}
The decision variable for (MINLP) is $\bchi = \{ \by, \bz \} \in \calY_a
 \times  \calZ_a$, where $\calY_a$ is a continuous variable set and
 $\calZ_a $ is a set of integers.
 It was shown in \cite{ruan-gao-minl} that for any given integer set $\calZ_a$, there exists a linear transformation
 $\bD_z :\calZ_a \rightarrow \mathbb{Z}  = \{ \pm 1\}^n$.
Thus, based on the unified  model (\ref{eq-unimod}), a general MINLP problem can be proposed as
\eb
(\calP_{mi}): \;\; \min \{ \Pi(\by,\bz) = \WW(\bD_y \by, \bD_z \bz) -  \la  \by ,  \bs \ra
 -   \la \bz ,  \bt \ra  \;  | \;\;(\by,\bz) \in \calY_c \times {\calZ}_c  \},
\ee
where
$\bff = ( \bs, \bt )$ is a given input,
$\bD\bchi = ( \bD_y \by, \; \bD_z \bz) \in \calW_y \times \mathbb{Z} $ is a multi-scale operator, and
\[
\calY_c = \{ \by \in \calY_a  \; | \;\; \bD_y \by \in \calW_y \},
\; \; \calZ_c = \{ \bz \in \calZ_a   | \; \bD_z \bz \in \mathbb{Z} \}.
\]
In $\calY_a$ certain linear constraints are given.  Since
 the set $\mathbb{Z}$ is bounded, by Assumption 1
 either $\WW:\calW_y \rightarrow \real$ is coercive or $\calW_y$ is bounded.
This general problem $(\calP_{mi}) $  covers many real-world applications, including  the so-called fixed cost problem \cite{g-r-s-fix}.
It must be emphasized that the integer constraint $\bw_z = \bD_z \bz \in \mathbb{Z}$ is a constitutive condition 
governed by the physical property of the system, it must  be relaxed by the canonicality. 
Let 
\eb
\beps =  \bLam_z(\bz) = (\bD_z \bz)  \circ  (\bD_z \bz)   \in \calE_z = \real^n_+,
\ee
where $\bx \circ \by = \{ x_i y_i\}^n$ is the Hadamard product in $\real^n$, thus 
the integer constraint in $\mathbb{Z}$ can be relaxed naturally by the
canonical function  $\Psi(\beps) = \{ 0 \mbox{ if } \beps \le  \be , \; \infty \mbox{ otherwise} \}$,
where $\be = \{1\}^n$ \cite{gao-jimo07}. 
Therefore, the canonical form of $(\calP_{mi})$ is
\eb
\min \{ \Pi(\by,\bz) = \Phi(\bLam(\by,\bz)) + \Psi(\bLam_z(\bz)) - \la  \by ,  \bs \ra  -  \la \bz,  \bt \ra
 \;  | \;
 \by\in \calY_c  \}.
\ee
Since the  canonical function $\Psi(\beps)$  is convex, semi-continuous,
its Fenchel conjugate is
\[
\Psi^\sharp(\bsig) = \sup\{ \la \beps ; \bsig \ra - \Psi(\beps) | \beps \in \real^n \}
=\{ \la \be ; \bsig \ra  \;\; \mbox{ if } \bsig \ge 0, \;\; \infty \mbox{ otherwise} \}.
\]
The generalized canonical duality relations (\ref{eq-gcdr}) are
\eb
 \bsig \ge 0 \;\; \; \Leftrightarrow \;\;\;  \beps  \le  \be   \;\; \Leftrightarrow  \;\;\;    \la  \beps  - \be ; \bsig \ra = 0   .
\ee
The complementarity shows that  the canonical integer constraint $\beps = \be $ can be naturally
relaxed by the
   $\bsig >  {\bf 0}$ in continuous space.
   Thus, if $ \bxi = \bLam(\bchi)$ is a quadratic homogenous operator and the canonical function
   $\Phi(\bxi)$ is convex on $\calE_a$,
    the canonical dual to $(\calP_{mi})$ is
   \eb
   (\calP^d_{mi}): \;\; \max \left\{ \Pi^d(\bvsig, \bsig) = - \half \la [\bG(\bvsig, \bsig)]^+  \bff , \bff \ra
   - \Phi^\sharp(\bvsig) - \la \be ; \bsig \ra | \;\; (\bvsig , \bsig) \in \calS^+_c \right\},
   \ee
   where $\bG(\bvsig, \bsig)$ depends on the quadratic operators  $\bLam(\bchi)$ and $\bLam_z(\bz)$,
   $\calS^+_c$ is a convex  open set
       \eb
       \calS^+_c = \{ (\bvsig, \bsig) \in \calE^*_a \times \real^n_+ | \;\; \bG(\bvsig, \bsig) \succeq 0, \;\; \bsig > 0 \}.
       \ee

   The canonical duality-triality theory has be used successfully for solving mixed integer programming problems
\cite{gao-ruan-jogo10,g-r-s-fix}. Particularly, for  the quadratic integer programming problem
\eb
(\calP_{qi}): \;\;\; \min \left\{ \Pi(\bx) = \half \bx^T \bQ \bx - \bx^T \bff | \;\; \bx \in \{ -1 , 1 \}^n \right\},
\ee
we have $\calS^+_c = \{ \bsig\in \real^n_+ | \;\; \bG(\bsig) = \bQ + 2 \Diag(\bsig) \succeq 0, \;\; \bsig > 0 \}$  and
\eb
(\calP^d_{qi}): \;\;\;\max \left \{ \Pi^d(\bsig) = -\half   \bff^T [\bG(\bsig)]^{+}  \bff  - \be^T  \bsig    | \;\; \bsig \in \calS_c^+ \right\}
\ee
which can be solved easily if   int$\calS_c^+ \neq \emptyset $.
Otherwise,  $(\calP_{qi})$  could be NP-hard since $\calS^+_c$ is an open set, which is a conjecture proposed in \cite{gao-jimo07}.
In this case,  $(\calP_{qi})$ is canonically dual to
  an unconstrained nonsmooth/nonconvex  minimization problem
\cite{gao-cace}.

\subsection{Relation with  SDP Programming}
Now let us show  the relation between the canonical duality theory and the popular semi-definite programming
relaxation.
\begin{thm} Suppose that $\Phi:\calE_s \rightarrow \real$ is convex and  $\barbvsig \in \calE^*_a$ is a solution of the problem
\eb
(\calP^{sd}): \;\;\; \min  \{ g + \Phi^*(\bvsig) \} \;\; s.t. \left( \begin{array}{cc}
\bG(\bvsig) & \bff \\
\bff^T & 2 g \end{array} \right)
\succeq 0 \;\; \forall \bvsig \in \calE^*_a , \;\; g \in \real,
\ee
 then
$\bchi = [\bG(\bvsig)]^+ \bff$ is a global minimum solution to the nonconvex problem $(\calP)$.
\end{thm}
{\bf Proof}.
The problem  $(\calP^d)$ can be equivalently written in the following problem (see \cite{z-g-y})
\eb
\min \left\{ g + \Phi^*(\bvsig) |\;\; g \ge  \Gap^\Lam(\bvsig), \;\; \bG(\bvsig) \succeq 0 \;\; \forall \bvsig \in \calE^*_a\right \}.
\ee
Then, by using the Schur complement Lemma \cite{schur}, this problem is equivalent to  $(\calP^{sd})$.
 The theorem is proved by the triality theory. \hfill $\Box$

It was proved   \cite{gao-ruan-jogo10} that for the same problem $(\calP_{qi})$, if we use different geometrical operator
 \[
 \Lam(\bx) =  \bx \bx^T \in \calE_a = \{ \bxi \in \real^{n\times n}| \; \bxi = \bxi^T, \;\; \bxi \succeq 0, \;\; \rank \bxi = 1, \;\; \xi_{ii} = 1 \;\; \forall i = 1, \dots, n\},
 \]
 and the associated canonical function
$ \Phi(\bxi) = \half  \la  \bxi ; \bQ \ra +  \{ 0  \mbox{ if } \; \bxi \in \calE_a ,
 + \infty  \mbox{ otherwise} \}$,
where $\la \bxi ; \bvsig \ra = \tr (\bxi^T \bvsig)$, we should obtain the same canonical dual problem $(\calP^d_{qi})$.
Particularly, if $\bff = 0$, then ($\calP_{qi}$)  is  a typical linear semi-definite programming
\[
\min \half \la  \bxi ; \bQ \ra \;\; s.t. \; \bxi \in \calE_a .
\]
Since $\calE_a$ is not bounded and there is no input, this problem is not properly posed, which could have
 either no solution or multiple solutions for a given  indefinite $\bQ = \bQ^T$.

\subsection{Relation to Reformulation-Linearization/Convexification Technique}
The {\em Reformulation-Linearization/Convexification Technique} (RLT) proposed  by H. Sherali and C.H. Tuncbilek \cite{han-tun92}
 is one well-known novel approach for efficiently solving  general  polynomial programming problems.
 The key idea of this technique is also to introduce a geometrically nonlinear operator $\bxi = \Lam(\bx)$
 such that the higher-order polynomial object $\WW(\bx)$ can be reduced to a lower-order polynomial $\Phi(\bxi)$.
 Particularly, for the quadratic minimization problems with linear inequality constraints in $\calX_a$:
\eb
 (\calP_{q}): \;\;\min \left\{  \Pi(\bx) =    \half \bx^T \bQ \bx - \bx^T \bff |\;\; \bx \in \calX_a \right\},
 \ee
by choosing  the   quadratic transformation
 \eb
\bxi = \Lam(\bx) =  \bx \overrightarrow{\otimes} \bx   \in  \calE_a  \subseteq  \real^{n\times n} ,
\;\; i.e., \;\; \bxi = \{ \xi_{ij} \} = \{ x_i x_j \}, \;\; \forall  1 \le i \le j \le n,
\ee
where $\overrightarrow{\otimes} $ represents the Kronecker product (avoiding symmetric terms, i.e. $\xi_{ij} = \xi_{ji}$),
   the quadratic object $\WW(\bw)$ can be reformulated as the following {\em first-level RLT linear relaxation}:
\eb\label{wphi}
\WW(  \bx) = \half \bx^T \bQ \bx =  \half \sum_{k=1}^n q_{kk} \xi_{kk} +   \sum_{k=1}^{n-1} \sum_{l=k+1}^n q_{kl} \xi_{kl} = \Phi(\bxi).
\ee
The linear   $\Phi(\bxi)$  can be considered as a special
  canonical function  since  $\bvsig = \partial \Phi(\bxi)$ is  a constant and
  $\Phi^*(\bvsig) = \la \bxi ; \bvsig \ra - \Phi(\bxi) \equiv 0$ is uniquely defined.
Thus,  using $\Phi(\bxi) = \la \bxi ; \bvsig \ra $ to replace $\WW( \bx)$
and considering  $\bxi$ as an independent variable, the problem $(\calP_{q})$ can be relaxed by
the following {\em RLT linear program}
\eb
(\calP_{RLT}): \;\; \min \left\{   \Phi(\bxi) - \la \bx , \bff \ra | \;\; \bx \in \calX_a, \;\; \bxi \in \calE_a \right\}.
\ee
Based on this RLT linear program, a branch and bound algorithm
was designed  \cite{han-tun}.
 It is proved that  if $(\barbx, \barbxi)$  solves $(\calP_{RLT})$,
 then its objective value yields a lower bound of $(\calP_q)$
  and $\barbx$ provides an upper bound
for $(\calP_q)$. Moreover, if  $\barbxi  = \Lam(\barbx) = \barbx \overrightarrow{\otimes} \barbx$,
then $\barbx$ solves $(\calP_q)$.

  This technique has been significantly adapted
along with supporting approximation procedures to solve a variety of more general nonconvex constrained
optimization problems having polynomial or more general factorable objective and constraint functions \cite{sherali}.

By the fact that for any symmetric $\bQ$, there exists $\bD \in \real^{n\times m}$ such that $\bQ = \bD^T \bH \bD$   with
$\bH = \{ h_{kk} = \pm 1, \;\; h_{kl} = 0 \;\; \forall k\neq l\} \in \real^{m\times m}$,
the canonicality condition (\ref{wphi}) can be simplified as
\eb
\WW(\bD \bx) = \half (\bD\bx)^T \bH (\bD\bx) =  \half \sum_{k=1}^m h_{kk} \xi_{kk} = \Phi(\bxi), \;\; \bxi =\Lam(\bx) = (\bD \bx) \overrightarrow{\otimes} (\bD \bx) \in \real^{m\times m}.
\ee
Clearly, if the scale $m\ll n$, the problem $(\calP_{RLT})$ will be much easier than the problems using the
geometrically nonlinear operator $\bxi =  \bx \overrightarrow{\otimes} \bx$.
Moreover, if we using the Lagrange multiplier $\bvsig \in \calE^*_a = \{ \bvsig \in \real^{m\times m}
| \;\; \la \Lam(\bx) ; \bvsig \ra \ge 0 \;\;\forall \bx \in \real^{n} \}$
to relax the ignored geometrical condition
$\bxi = \Lam(\bx)$ in $(\calP_{RLT})$, the problem $(\calP_q)$ can  be equivalently relaxed as
\eb
(\calP_\Upsilon): \;\; \min_{\bx, \bxi} \max_{\bvsig} \left \{  \Upsilon(\bx, \bxi, \bvsig) = \Phi(\bxi) + \la  \Lam(\bx) - \bxi ; \bvsig \ra - \la \bx ,\bff \ra | \;\; \bx \in \calX_a, \;\; \bxi \in \calE_a , \;\; \bvsig \in \calE_a^* \right\}.
\ee
Thus, if $(\barbx, \barbxi,\barbvsig)$ is a solution to $(\calP_\Upsilon)$, then $\barbx$ should be a solution to $(\calP_q)$.
By using the sequential canonical quadratic transformation $ \Lam(\bx) = \Lam_p (\dots (\Lam_1(\bx) \dots )$ (see Chapter 4, \cite{gao-dual00}), this technique can be used for solving general global optimization problems.

\section{Symmetry, NP-Hardness and Perturbation Methods}
\begin{figure}[h]
\begin{picture}(3,2)
\setlength{\unitlength}{.5cm}
\put(2,-10){{\large \psfig{file=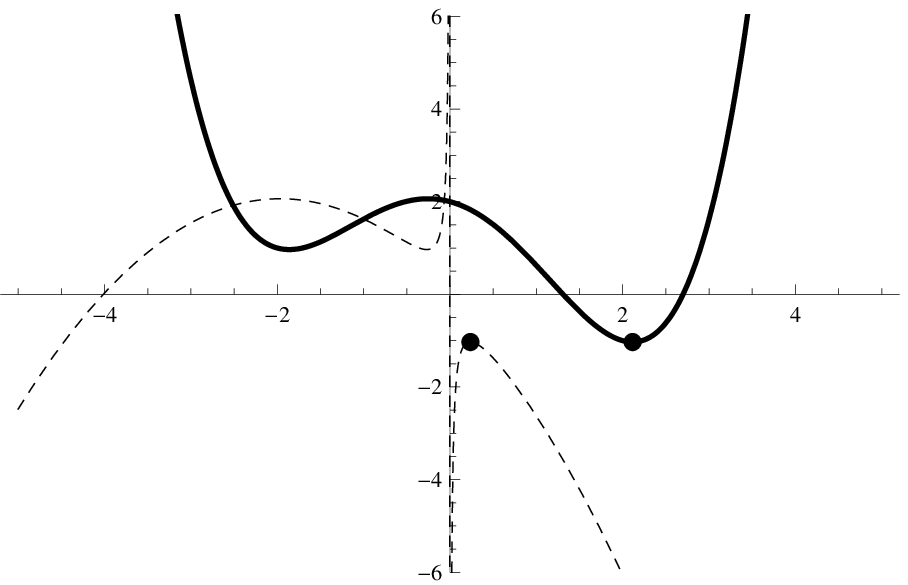,height=3.2cm, width =4cm}}}
\put(10,-10){{\large \psfig{file=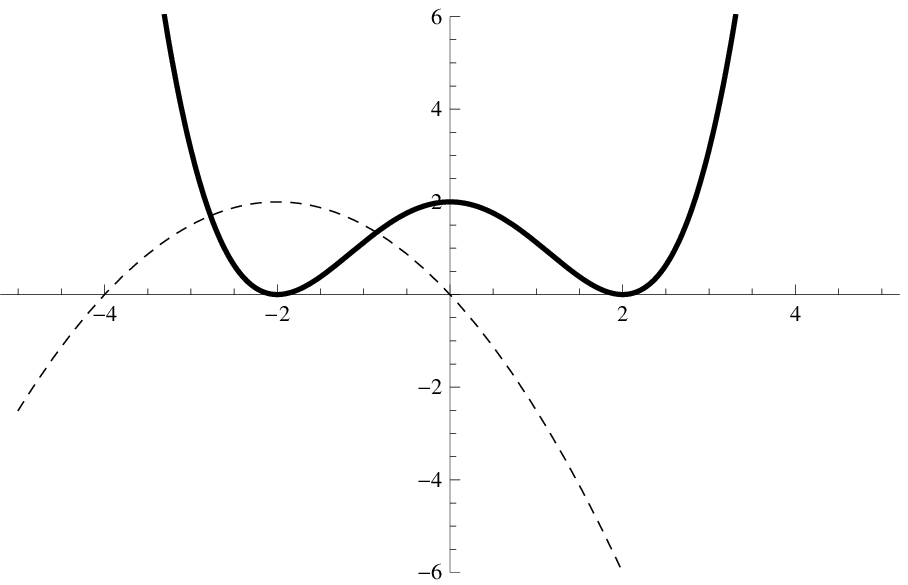,height=3.2cm, width =4cm}}}
\put(20,-10){{\large \psfig{file=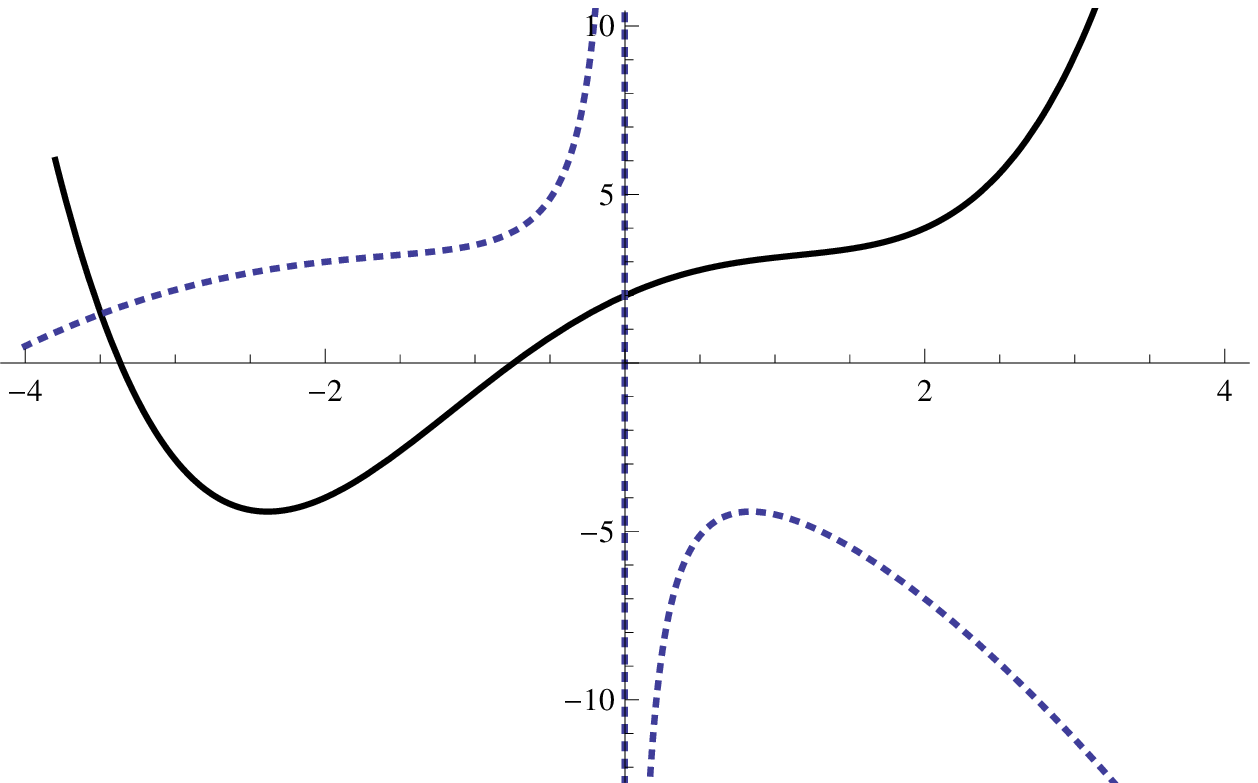,height=3.2cm, width =3cm}}}
\put(2,-9.5){(a) $f=0.5 $}
\put(10,-9.5){(b)   $f = 0 $}
\put(18,-9.5){(c)   $f = -2.0$}
\end{picture}\vspace{4.8cm}
\caption{ Graphs of   $ \Pi (\bx)$  (solid) and
   $\Pi^d(\vsig)$ (dashed) ( $\alpha= 1$, $\lam=2$ ) }\label{onedim}
\end{figure}

The concept of symmetry is closely related to the duality and, in certain sense, can be viewed as a
 {\em geometric duality}.
 Mathematically,  symmetry means invariance under  transformation.
 By  the  canonicality,
the object $\WW(\bw)$ possesses  naturally  certain symmetry.
If the subject  $\FF(\bchi) = 0$, then  $\Pi(\bchi) =\WW(\bD \bchi) =\Phi(\Lam(\bchi))$ and $(\calP)$
should have either a trivial solution  or multiple solutions due to the symmetry.
   In this case $\Pi^d(\bvsig) = - \Phi^*(\bvsig)$ is concave and,   by the triality theory,
   its critical   point $\barbvsig \in \calS^-_c$ is a global maximizer,
   $\barbchi = [\bG(\barbvsig)]^+ \bff = 0$ is the biggest local maximizer of $\Pi(\bchi)$,
   while the global minimizers must be $\barbchi(\barbvsig)$  for those  $\barbvsig \in \partial \calS_c^+$
   such that $\Pi^d(\barbvsig) = \min \{  - \Phi^*(\bvsig) | \;\; \det\bG(\bvsig) = 0 \;\; \forall \bvsig \in \calS_c\}$.
   Clearly,  this nonconvex constrained concave minimization problem could be really NP-hard.
Therefore, many well-known NP-hard problems in computer science and global optimization are not well-posed problems.
Such as the {\em max-cut problem},
which  is a special case of quadratic integer programming problem $(\calP_{qi})$.
 Due to the symmetry $\bQ=\bQ^T$ and  $\bff = 0$,
 its canonical dual problem has multiple solutions on the boundary of
 $\calS^+_c$. The problem is considered as NP-complete even if $Q_{ij} = 1$ for all edges. 
Strictly speaking, this is not a real-world problem but only a  geometrical model.
Without sufficient  geometrical constraints in $\calX_a$, the graph is not physically fixed
 and any rigid motion is possible.
   However, by adding a linear perturbation $\bff \neq 0$,
   this problem can be solved efficiently  by the canonical duality theory  \cite{wang-ea}.
   Also it was proved by the author \cite{gao-cace,gao-ruan-jogo10} that the general quadratic integer problem
   $(\calP_{qi})$ has a unique solution as long as the input $\bff \neq 0$ is big enough.
  These results show that the subjective function plays an essential role for symmetry breaking to leads a well-posed  problem.
To  explain the theory and understand the NP-hard problems, let us consider a simple problem  in
$\real^n$:
\begin{eqnarray}
\min \left\{  \Pi(\bx)=\half \alpha (\half\|\bx\|^2-\lam)^2-\bx^T \bff  \;\;\;   \forall \bx \in \real^n \right\} ,\label{examdoub}
\end{eqnarray}
where $\alp, \lam > 0$ are given parameters.
 Let   $ \Lam(\bx)=\half\|\bx\|^2  \in \real$, the canonical dual function is
\begin{eqnarray}
\Pi^d(\vsig)=   - \half \vsig^{-1}   \|\bff \|^2  -\lam \vsig  -\half \alpha^{-1} \vsig^2 ,
\end{eqnarray}
 which is defined on $\calS_c = \{ \vsig \in \real| \;\; \vsig \neq - \lam , \;\; \vsig = 0 \mbox{ iff } \bff = 0\}$.
 The criticality condition  $\partial \Pi^d(\vsig) = 0$ leads to a canonical dual equation
 \eb\label{eq-cde}
(\alpha^{-1} \vsig+\lam)\vsig^2=\half \|\bff \|^2 .
\ee
This cubic equation has at most three real solutions satisfying $\vsig_1 \ge 0 \ge \vsig_2\ge \vsig_3$,
and, correspondingly,    $\{\bx_i=\bff/\vsig_i \}$ are three critical points of $\Pi(\bx)$. By
the fact that $\vsig_ 1 \in \calS^+_a = \{ \vsig \in \real\; |\; \vsig \ge  0 \}$,
  $\bx_1$ is a global minimizer of $\PP(\bx)$. While for $\vsig_2, \vsig_3 \in \calS^-_a=   \{ \vsig \in \real\; |\; \vsig <  0 \}$, $\bx_2$ and $\bx_3$  are local min (for $n=1$) and local max of $\Pi(\bx)$,
   respectively (see   Fig. \ref{onedim}(a)).\vspace{-1.6cm}

If we let $\bff= 0$, the graph of $\Pi(\bx)$ is symmetric (i.e. the so-called double-well potential or the Mexican hat for $n=2$
\cite{gao-amma03})
 with  infinite number of global minimizers
 satisfying $\| \bx \|^2  = 2 \lam$.
In this case, the canonical dual $\Pi^d (\vsig) = - \half \alp^{-1} \vsig^2 - \lam \vsig$ is strictly concave
with only one critical point (local maximizer) $\vsig_3 = -  \alp \lam  < 0 $.
The corresponding solution $\bx_3 = \bff /\vsig_3 = 0$ is a local maximizer.
 By the canonical dual equation (\ref{eq-cde})
we have $\vsig_1 = \vsig_2 = 0$ located on the boundary of $\calS^+_a$,
which  corresponding to the two global minimizers   $x_{1,2} = \pm \sqrt{2 \lam}$ for $n=1$, see Fig. 1 (b).
If we let $f =-2$, then the graph of $\Pi(\bx)$ is quasi-convex with
only one critical point. In this case, (\ref{eq-cde}) has only one solution $\vsig_1 \in \calS^+_c$
(see Fig.1 (c)).

 \begin{Conj}
For any given properly posed problem  $(\calP)$ under the   Assumption 1,
  there exists a  constant $f_c > 0$ such that
 $(\calP^d)$   has a unique solution in $\calS^+_c$ as long as   $\|\bff\| \ge f_c$.
 \end{Conj}
 
This conjecture shows that any properly posed problems are  not NP-hard if the input  $\|\bff\|$ is big enough.
Generally speaking,   most  NP-hard problems have multiple solutions located either on
  the boundary  or the outside of $\calS^+_c$.
Therefore,  a quadratic perturbation method can be suggested as
\begin{eqnarray*}
\Xi_{\delta_k} (\bchi, \bvsig)  =  \Xi(\bchi, \bvsig) + \half \delta_k \| \bchi - \bchi_k \|^2
 =   \half \la \bchi , \bGd(\bvsig) \bchi \ra - \Phi^*(\bvsig) - \la \bchi, \bfd \ra + \half \delta_k \la \bchi_k, \bchi_k \ra ,
\end{eqnarray*}
where $\delta_k > 0, \;  \bchi_k \; (k = 1, 2, \dots) $ are  perturbation parameters, $\bGd(\bvsig) = \bG(\bvsig) + \delta_k \bI $, and $\bfd  = \bff+ \delta_k \bchi_k$.
Thus, the original canonical dual feasible space $\calS_c^+$ can be enlarged to
$\calS_{\delta_k}^+ = \{ \bvsig \in \calS_c | \; \bGd(\bvsig) \succ 0 \}$
such that  a perturbed canonical dual problem can be proposed as
\eb
(\calP^d_k): \;\;  \max \left\{ \min \{ \Xi_{\delta_k} (\bchi, \bvsig) |\;\; \bchi \in \calX_a \} | \; \; \bvsig \in \calS_{\delta_k}^+ \right\}.
\ee
Based on this problem,  a canonical primal-dual algorithm has been developed with successful applications
for solving sensor network optimization problems \cite{ruan-gao-pe} and chaotic dynamics \cite{lat-gao-chaos}.\vspace{-.2cm}

\section{ Challenges and Breakthrough}
 Now let us turn our attention to the most recent  challenges from C. \z,
 who claimed in  the beginning of his new paper \cite{z16}:

 ``We believe one of the main aims of [5] is to question our counter-examples in [7], [6],
[8], [9], [10], [11], [12], [15], [16], [13], [14]."

Readers can easily verify this false assertion  by checking the first version\footnote{posted at  http://arxiv.org/abs/1310.2014} of \cite{lat-gao-opl}.
Although the canonical duality-triality theory has been repeatedly  challenged by
\za with two co-workers in the 11 papers [6-14] and
the author has been invited by journal editors to write responses, he didn't do it as
   the  mistakes in these papers are so basic and should be easily understood by experts in the
  communities. Also  the author believes that the  {\em truth needs no defending}.
The only  half  paragraph of comments on these  false challenges added in the revision of
\cite{lat-gao-opl}  was based on the  reviewers' comments.
However,  even  these added comments are considered by \za as one of  main aims of \cite{lat-gao-opl}.
Indeed, \za  has been fully alarmed on  all author's recent papers
and anxious to against  any comments on his basic mistakes (see \cite{z14}).
Unfortunately,  without correct motivation
 and necessary knowledge  for understanding the canonical duality theory,
more  mistakes have been  produced  both ethically  and mathematically.\vspace{-.1cm}

\subsection{Mathematical Mistakes}
  \z's  first aim in \cite{z16} is to point out   ``a false
assertion, a not convincing proof, a non adequate application of a result of the paper, several
inconsistencies".

The so-called  ``false assertion" discovered by \za
 is only  a sentence in author's paper \cite{lat-gao-opl}, i.e.
{\em any convex quadratic function is objective}.
 By using a ``counterexample"
 $W(u,v) = u^2 + 2 v^2$, \za proved  that this assertion is false.
  However, it was written clearly in \cite{lat-gao-opl} (equation (2) on page 3)
   that the objective function
  must be in the form of $\WW(\DD\uu)$. If
   we simply let $D = \Diag (1, \sqrt{2})$ and $\bw  = D (u, v)^T = (u, \sqrt{2} v)^T$, then
    $W(\bw) = \|\bw \|^2 $ is truly  an objective function in $\real^2$.
    This basic mistake shows that  \za  works only on the one-scale artificial optimization problem (1) and
    dose not understand the
    multi-scale modeling of (\ref{eq-gs}),  a foundation not only for the canonical duality theory,
    but also for entitle applied mathematics
    beautifully presented  in G. Strang's textbook \cite{strang} from the first chapter of linear algebra to the last content on   optimization\footnote{MIT's online teaching project was started
    from Gil Strang's this textbook}.

    In the same section, \za complained that the
    Definition 6.1.2 of objective function  in \cite{gao-dual00}
    is not clear at least because $\Oo\times \calW_a$ does not seem
to be an ``objective set". This complain shows that \za is confused about the difference between
optimization in $\calX_a \subset \real^n$  and   nonlinear analysis in continuous space $\calW_a \subset \calL^p[\Oo; \real^d]$,
 where $\bx \in \Oo$
is not an unknown
and, as discussed in Section 2.3 of this paper,
the stored energy $\WW(\bw) = \int_\Oo \UU(\bx, \bw) \dO$ is objective
iff the stored energy density  $\UU(\bx, \bw)$ is an objective function of $\bw \in \calW_a$.
By any numerical approximation, $\bx$ will be disappeared such that  $\UU(\bx, \bw) \simeq \UU_k(\bw)$ and $\WW(\bw) = \sum_k \UU_k(\bw)$.

Regarding the so-called ``not convincing proof",
serious researcher    should provide  either a convincing proof or a disproof, rather than a complaint.
Note  that  the canonical dual variables $\bsig_0$ and $\bsig_1$ are
 in  two  different levers  (scales) with totally different physical units\footnote{Let us consider Example 1 in \cite{lat-gao-opl}. If the unit for $x$ is the meter $(m)$ and for $q$ is $Kg/m$, then
 the units for the Lagrange multiplier $\mu$ (dual to the constraint $g(x) = \half (\half x^2 - d)^2 - e$)
 should be $Kg/m^3$ and for $\sig$ (canonical dual to $\Lam(x) = \half x^2$) should be $Kg/m$, respectively, so that
 each terms in $\Xi_1(x,\mu,\sig)$ make physical sense.},
it is completely  wrong to consider $(\bsig_0, \bsig_1)$ as one  vector and  to discuss the
  concavity of $\Xi_1(x, (\cdot,\cdot))$ on  $\calS^+_a$.
 The condition  ``$\calS^+_a$ is convex" in Theorem 2 \cite{lat-gao-opl}
should be understood in the way that  $\calS^+_a$ is convex in $\bsig_0$ and $\bsig_1$,
respectively, as emphasized  in Remark 1 \cite{lat-gao-opl}.
Thus, the proof of Theorem 2 given in \cite{lat-gao-opl} is indeed convincing
by simply using the classical saddle min-max duality for $(\bx, \bsig_0)$ and $(\bx,\bsig_1)$, respectively.

Regarding the non adequate application, it is easy to check
$x_1$ is  a global minimizer since $(\mu_1, \sig_1) \in \calS^+_a$ and
$\calS_a^+ = \{ (\mu,\sig) \in \real \times \real| \;\; \mu > 0, \; q + \mu \sig > 0 \}$
is convex in $\mu$ and $\sig$ for any given $q \in \real$.
Thus,  Example 1  is indeed an adequate application of Theorem 2 in \cite{lat-gao-opl}.
Dually, \z's modified stored energy function $W(\DD x) = - \half x^2$   in \cite{z16}
is totally artificial (against the basic law in current physics and
Assumption A1.1,  A1.3), it is not surprise  to have int$\calS^+_a = \emptyset$, which is the case for many artificially produced  NP-hard problems.

The several inconsistencies discovered by \za proved the triality in the mathematical modeling $(\calP_0)$, i.e.
 the subjective function $\FF(\uu) = - \la \uu, \barbu^* \ra $   is necessary for any real-world problems,
 including the  paper writing.
 In physics, the input $\barbu^*$ could be defects,  white  noise, or random dislocation and charges,  etc (see \cite{g-l-v}).
 Although there are many such  ``inconsistencies"  and even mistakes in his publications  including \cite{gao-dual00,gao-strang89},
the author never write any erratum as he has been busy in search natural  beauty but  realized
 that nothing  is perfect in this real world.
 Indeed, readers can easily find many typos and inconsistencies in \z's papers,
 especially in his  current short paper \cite{z16}. \vspace{-.3cm}

\subsection{Non-Mathematical  Mistakes}
 \z's second aim in \cite{z16} is quibbling about  the basic mistakes in their false challenges  published in [6-14].
 The author believes that by  now  serious researchers should have a clear understanding on the canonical duality-triality  theory,
 its  history and unified applications to multidisciplinary fields,
 so he  has no intention to argue with \za on  tedious  issues
   but to point out some hidden  truths
in the arguments and ethical mistakes made in  \cite{z16}.

a*). In Section 2 a), \za wrote:  ``it is not possible to find a word containing `bjective' in Gao-Strang's paper [2].
So we have not any reason to guess that $W$ has to be `objective' and $F$ to be `subjective'.
 In fact D.Y. Gao introduced 11 years later the notion "objective function" in
 his book [1], and also refers to a book by P.G. Ciarlet published in 2013". 

Objectivity is a basic concept in continuum physics and is usually discussed in the beginning of textbooks,
say page 8 in \cite{mar-hug} (published in 1983). 
Any one who ever took a graduate course of  nonlinear field theory should know that
the stored energy $\WW$  must be objective, therefore, most research papers don't mention this basic requirement except some review/survey articles on interdisciplinary topics, say \cite{ball}.
P.G. Ciarlet's new  book (2013) on nonlinear analysis \cite{ciarlet} is
 based on his well-known book on nonlinear elasticity published in 1988,
 in which, this basic requirement is  called the {\em  axiom of objectivity} (see page 101 \cite{ciarlet88}).
 The mathematical definition of the objectivity given in author's  book
  (Definition 6.1.2 \cite{gao-dual00}) is
 based on  Axiom 3.3-1 in \cite{ciarlet88}, which is correct and clear for any mathematicians who
  know    the difference between the coordinates $t, \bx \in \Oo$ and the Lagrangian coordinates  $\bu(t, \bx)$.
 In Gao-Strang's paper, the geometrical constraint in $\calU_a$ is relaxed by the indicator $\Psi_{\calU_a}(\bu) = \{ 0 \mbox{ if } \bu \in \calU_a, \; \infty \mbox{ if } \bu \notin \calU_a\}$
 such that the ``subjective function" is written as the {\em external super-potential}
 $\FF(\bu) = -\la \bu, \barbu^* \ra + \Psi_{\calU_a}(\bu)$ (see equation (85) in \cite{gao-strang89}).
 It is a common  sense that the minimum potential energy principle holds only for
 the so-called {\em dead loading} systems, i.e. the input $\barbu^* = \partial \FF(\bu)$ is independent of the output $\bu$.   Thus the  external energy $\FF(\bu)$ must be linear on its effective domain.
 Otherwise, the system is not conservative and the traditional  variational method can't be applied.
 In order to solve such  problems, a so-called {\em rate-variational method} was proposed by Gao and Onat in 1990 \cite{gao-onat}.
 However, 
  M.D. Voisei and C. \za oppositely
choose   piecewise linear function $W$  and quadratic function  $F$
as counterexamples to against Gao and Strang's paper
 with six conclusions including [10]
 ``About the (complementary) gap function one can conclude that it is useless
 at least in the current context. The hope for reading an optimization theory with diverse applications is ruined \dots"
 Clearly, this is a duality mistake. 

 \za continues in a) with: ``we didn't find any mention about `objective function' with the meaning from [1, Definition 6.1.2] until 2010 in Gao's articles, that is before submitting all our papers on Gao's works".
However, in his open letter to the author\footnote{\za (2012): ``Open letter to David Yang  Gao" (version 1) posted on his university web page},
 \za wrote: ``In 2006 I proposed my former student R. Strugariu to study your theory. \dots, Z2G 27 May 2008: Some time ago I bought your book mentioned below and I began to read it".
Clearly,  \za indeed  read the  book \cite{gao-dual00} before 2010, where the objective function
 is mathematically defined.
 This contradiction proves that \za does not  tell the truth in the arguments.

 b*). \za wrote in b): ``Indeed, we proved that practically all statements called `triality theorem' in Gao's papers
published before 2010 are false".

As discussed in Section 3.3, the triality theorem was correctly proposed in 1997 from a
 post-buckling problem, where $\dim \Pi = \dim \Pi^d$ \cite{gao-amr97}.
Mathematical proof  was given in 2000  \cite{gao-na00} for nonconvex variational problems.
Generalization to global optimization was made in 2000. 
Sooner in 2002 \cite{gao-opt03, gao-amma03} the author  discovered  that  the double-min duality does not hold
in its strong form if $\dim \Pi \neq  \dim \Pi^d$. But the triality theorem was still presented correctly in the
``either-or" form since the double-max duality is always true.
Careful readers can find immediately that instead of the original paper \cite{gao-na00}, its applications
\cite{gao-ogden-zamp,gao-ogden-qjmam} were challenged
 and the three key papers \cite{gao-amr97,gao-opt03, gao-amma03}
were never cited by \za and his co-workers in [6-14].

It is now necessary to exam how they proved  the ``false" of the triality theorem.

Among the 11 papers [6-14], seven of them against the triality theorem in global optimization.
The so-called proof  is mainly  a set of ``counterexamples"  in [7,8,11,14-16] with  $\dim \Pi \neq  \dim \Pi^d$,
which  is exactly the ``additional condition"
discovered by the author in  \cite{gao-opt03, gao-amma03}.
This should be the only reason why  \cite{gao-opt03, gao-amma03} were not cited by these people.
One counterexample in [9] is  the case
that  $\barbvsig \in \partial \calS^+_c $, so Voisei and  \za concluded: ``The consideration of the function
$\Xi$ is useless, at least for the problem
studied in \cite{gao-yang}"\footnote{The reference \cite{gao-yang} is [3] in [9]}.
However, it was proved in \cite{mora-gao-surface} that by simply using a line perturbation $\bff \neq 0$, this
so-called ``counterexample" can be solved nicely by the triality theorem to obtain all global optimal solutions.

The rest papers [6,10,12,13] challenged the triality theory in nonconvex analysis.
As discussed above, ``counterexamples"  in [6] are simply using linear $W$ and nonlinear $F$ to
 against Gao-Strang's paper \cite{gao-strang89}.
It is very interesting to point out that, instead of arguing the mathematical proof given in \cite{gao-na00},
Voisei and  \za  challenged  Gao and Ogden's papers \cite{gao-ogden-zamp,gao-ogden-qjmam}, which are applications of \cite{gao-na00}.
 By using   ``a thorough analysis" and a convention $0/0 :=0$ in measure theory,
they proved in [12,13]
 that the main result in triality theory is false and if $\beta \neq 0$ (i.e.   $\bff \neq 0$)
 the canonical dual $\Pi^d(\bvsig)$ is not well-defined.
 Unfortunately, they don't know a basic truth guaranteed by  the canonical duality theory, i.e.
 if  $\bff \neq 0$, then $\bvsig \rightarrow 0$  can never happen (see Equation (\ref{eq-cde})),
 otherwise, Newton's third law will be violated.
 Therefore, $\Pi^d(\bvsig)$ is indeed well-defined.
 Also, the constitutive law $\bvsig = \partial \Phi(\bxi)$
for phase transitions in continuum mechanics   holds only  at the scale about $\bvsig \sim 10^{-3}m$.
Even the quantum mechanics can't reach the scale of zero measure. Therefore,  their  papers  [12,13] were rejected.
However,  Voisei and  \za couldn't understand  their  basic mistakes  and had serious arguments with
editors of these two decent  mathematical physics journals\footnote{see comments and links posted on  http://arxiv.org/abs/1202.3515
and http://arxiv.org/abs/1101.3534}.

The most funny  mistake  ever made by    \za and his co-workers
could be the one in their  paper [6] published in a   dynamical  systems journal.
As it is known that the bi-duality was first proposed and proved
by the author for  convex Hamiltonian systems \cite{gao-dual00}, 
where the Lagrangian must be in its standard form $L(\bq, \bp)$, i.e. Equation (\ref{eq-lag}) in Lagrangian mechanics.
Instead of finding any possible mistakes in author's proof,
 Strugariu, Voisei and  \za  created an
 artificial ``Lagrangian":  
\[
L(x, y) := -\half \alp \| x\|^2 - \half \beta \| y\|^2  + \la a , x \ra \la b, y\ra, \;\; \;\; (\mbox{Equation (1) in \cite{svz}})
\]
By using this  ``Lagrangian" as well as   the associated  ``total action" and its Legendre   dual
\begin{eqnarray*}
f(x) &=& \max \{ L(x,y) | \; y \in Y \} = -\half \alp \| x\|^2 + \half \beta^{-1} \la a, x \ra^2 \|b\|^2 \;\; \forall x \in X \\
 g(y) & =&   \max \{ L(x,y) | \; x \in X \} = - \half \beta \| y\|^2 + \half \alp^{-1} \la b, y\ra^2 \|a \|^2 \;\;\; \forall y \in Y
 \end{eqnarray*}
 they   produced a series of very strange counterexamples to against the  bi-duality theory in
  convex Hamiltonian systems   and  the triality theory in  geometrically nonlinear systems
 presented  respectively by the author in Chapters 2 and 3  \cite{gao-dual00}.
They claimed :``Because our counter-examples are very simple, using quadratic functions defined
on whole Hilbert (even finite dimensional) spaces, it is difficult to reinforce the
hypotheses of the above mentioned results in order to keep the same conclusions
and not obtain trivialities."

Clearly, the  quadratic function $L(x,y)$ is totally irrelevant to the Lagrangian
 $L(\bq,\bp)$ in Lagrangian mechanics and Gao's book \cite{gao-dual00}.
Without the differential operator $\DD = \partial_t$, the quadratic d.c.
function $f(x)$ (or $g(y)$)
 is defined on one-scale space $X$ (or $Y$) and is unbounded.
 Therefore, it's critical point  does not produce any motion.
 This basic mistake shows that these people don't have
  basic knowledge not only in  Lagrangian mechanics (vibration produced by the duality between the kinetic energy $\TT(\partial_t \bu)$ and the potential energy $\UU(\bu)$), but also in
  d.c. programming (unconstrained quadratic d.c. programming does not make any sense \cite{jin-gao}).
  It also shows that these people even don't  know what the Lagrangian coordinate is,
otherwise, they  would never  use a time-independent vector $x\in \real^n$ as an unknown in dynamical systems and
\za  wouldn't complain   the definition of the objectivity given in \cite{gao-dual00}.

Moreover, the triality theory  was developed from geometrically nonlinear systems,
where the geometrical operator $\Lam(\bu)$ must be nonlinear in order to have
canonicality condition (A1.2) and the triality theory (see \cite{gao-jogo00}).
By the fact that only  the geometrical nonlinearity can produce multiple local minimizers,
this is the reason why this terminology was emphasized in the title  of Gao-Strang's paper
 \cite{gao-strang89}.
 However, in \cite{svz} Strugariu, Voisei and  \za  choose either a null $\Lam(u) = 0$ (Example 3.4) or a linear
  $\Lam(u) = \la a, u \ra b$ (Example 3.5) as counterexamples to prove the false of the triality.
   These mistakes show  that these people really don't understand both the geometrical nonlinearity and the triality theorem.

  Even more, since   there is neither  input in $L(x,y)$ nor initial/boundary conditions in $X$,
  all counterexamples they produced are simply not problems but only artificial ``models".
  Since they don't follow  the basic roles  in mathematical modeling, such as the  objectivity,
   symmetry, conservation and constitutive  laws,
 etc,
  these artificial ``models" are very strange and even ugly (see Examples 3.3, 4.2, 4.4 \cite{svz}).

 All these conceptual  mistakes show  that \za and his two co-workers    don't know what they are doing:
without understanding the title (geometrical nonlinearity) of \cite{gao-strang89} and the basic contents (objectivity, stored energy and external energy) in nonlinear analysis, they published the paper [10] in {\em Applicable Analysis}
to  against Gao-Strang's work in nonlinear analysis;
 without necessary knowledge of Lagrangian mechanics they published
the paper \cite{svz}  in {\em Discrete \& Continuous Dynamical Systems-A}
to challenge  Gao's book on convex Hamiltonian systems.
Readers are suggested to check the special issues of these two journals  to understand why these  papers can be published.

 c*). \za wrote in c): ``I ask the authors of [5] to give precise
references where our counter-examples can be found in Gao's works (or elsewhere); otherwise
the statement ``these so-called counterexamples are not new, which were first discovered by
Gao" is a calumny."

All these counter-examples are simply using   the condition $\dim \Pi \neq \dim \Pi^d$
 to against  the double-min duality.
Indeed, such a type of counter-examples is too  simple, i.e. the double-well problem  (\ref{examdoub})
with $n \ge 2$, which  can be found easily  in
 author's book \cite{gao-dual00} and many articles, say \cite{gao-na00,gao-jogo00,gao-opt03,gao-jogo06}.
Precisely,  Example 5.1 in \cite{gao-opt03}  and  Example (2.14) in \cite{gao-amma03} are the counter-examples first discovered by
the author in 2003.
So it was written clearly in  Remark 1 on page 481 \cite{gao-opt03} and   Remark of Theorem 3  on
page 288 \cite{gao-amma03}
that the double-min duality (\ref{dmin}) holds ``under certain additional conditions".
It has been proved either in author's book \cite{gao-dual00} or in recent papers \cite{gao-wu11,gao-wu,chen-gao-jogo,gao-wu-jimo}
that $\dim \Pi = \dim \Pi^d$  is the only condition for the double-min duality.
Anyone who knows the logic will surely understand the counterexamples discovered by Gao in 2003 must be the same type as those listed in [6-14] by \za {\em et al}.
 The only reason  why the author didn't write down specifically the condition
 $\dim \Pi \neq \dim \Pi^d$ in \cite{gao-opt03,gao-amma03}
is that  he  was not sure if there is any other conditions, so  he was prefer to leave this  uncertainty as an open problem to readers, which is  author's philosophy as the old saying:``hidden harmony is stronger than the explicit one".
Serious researchers may  ask  why such  simple duplicated ``counter-examples" can be published
 repeatedly in  the  international journals without citing \cite{gao-opt03,gao-amma03}? 

As \za  indicated in his open letter, the author is indeed one of three reviewers for his paper \cite{11}
and   Gao's papers \cite{gao-opt03,gao-amma03} were pointed out in all the three reviewers' reports\footnote{The author thanks  these two reviewers for forwarding their reports. The second reviewer indicated specifically the Remark 1 in Gao's paper  \cite{gao-opt03} and he recommended these people to  ``present their `counter-examples' more as `examples ...' ".
The third reviewer wrote clearly ``but counterexamples had already been discovered by Gao in his earlier paper".}.
Unfortunately,  Voice and \za still refuse to cite these two key papers in their revision  \cite{11}
but simply deleted the similar sentence ``a correction of this theory is impossible without falling into trivia" as they conclude in [6].
   Since 2012  the author and his co-workers proved that
even if $\dim \Pi \neq  \dim \Pi^d $, the double-min duality still holds weakly in a beautiful symmetrical form \cite{chen-gao-jogo,gao-wu-jimo,mora-gao-memo}. \za knows these progress,
 at least the paper  \cite{gao-wu},
 so    his own   statement in  \cite{z16}  ``Indeed, we proved that practically all statements called `triality theorem' in Gao's papers
published before 2010 are false" is truly a calumny.

d*). In d) \za wrote: ``Indeed, we never cited Gao's papers [6,7]\footnote{i.e. Gao's paper \cite{gao-opt03,gao-amma03}}. The simple reason is that we learned about
the so called open problem from Gao's paper [6] (see footnote 4) from 2 reports on our paper
[11], received on 06.05.2010; at that moment (06.05.2010) all the 11 papers were already
submitted."

As we know that \za has begun interesting in the canonical duality theory at least from 2006. At that time, the author published a very few papers  in optimization journals.
It is difficult to believe that \za didn't read \cite{gao-opt03} before to criticize this theory.
Indeed, any people, if they simply check   \cite{14},  should know  immediately that
this paper by   Voisei and \za was submitted to the journal  on 27.4.2011, i.e. almost one year after ``that moment (06.05.2010)"
(also the author is a reviewer, both \cite{gao-opt03, gao-amma03} were mentioned in the   report),
but  neither \cite{gao-opt03} nor \cite{gao-amma03} was  cited by these people.
 This contradiction shows again that  \za does not tell the truth in \cite{z16}.
 So there is no need to continue this discussion.



All the conceptual and mathematical
 mistakes in this set of published/rejected  papers [6-14] by \za and his two co-workers
 show a significant gap between their ``thorough  mathematics"  and
 the {\em applicable mathematics} that
 the canonical duality-triality is based on.
 As V.I. Arnold concluded in his address \cite{arnold}:``A teacher of mathematics, who has not got to grips with at least some of the volumes of the course by Landau and Lifshitz, will then become a relict like the one nowadays who does not know the difference between an open and a closed set." \vspace{-.1cm}


\section{Conclusions}
Based on   necessary conditions and basic laws in physics,
a unified   multi-scale global optimization problem  is proposed in the canonical form: \vspace{-.1cm}
\eb
\Pi(\bchi) = \WW(\bD \bchi) + \FF(\bchi) = \Phi(\Lam(\bchi)) - \la \bchi, \bff \ra .\vspace{-.1cm}
\ee
The object $\WW$ depends only on the model and   $\WW(\ww) \ge 0 \;\forall \ww \in \calW_a$ is necessary;
 $\WW$ should be an objective function for physical systems, but it is not necessary for artificial systems (such as management/manufacturing processes and numerical simulations, etc).
  The subject $\FF$ depends on each properly posed problem and must satisfy $\FF(\bchi) \le  0$ together with  necessary  geometrical constraints for the output $\bchi \in \calX_a$ and
 equilibrium conditions for the input  $\bff  \in \calX^*_a$.
 The geometrical nonlinearity of $\Lam(\bchi)$ is necessary for nonconvexity in global optimization, bifurcation in nonlinear analysis,  chaos in dynamics, and NP-hardness in computer science.

Developed from large deformation nonconvex analysis/mechanics,
 the canonical duality-triality  is a precise mathematical theory with solid foundation in physics
 and natural root in philosophy,
 so it is naturally related to the traditional theories and powerful methods in
 global optimization and nonlinear analysis.
 By the fact that the canonical duality is a universal law of nature, this theory
 can be used not only  to  model real-world problems, but also for
solving a wide class of challenging problems in multi-scale complex systems.
The conjectures proposed in this paper can be used
  for understanding and clarifying NP-hard problems.

Both the linear operator $\bD $ and the geometrical admissible operator $\Lam$
can be generalized to the composition forms \vspace{-.2cm}
\eb
\bD = \bD_n \circ \cdots \circ \bD_1, \;\; \;\;\;  \Lam(\bchi) = \Lam_m(\Lam_{m-1} (\cdots (\Lam_1(\bchi ) ) \cdots ) )\vspace{-.2cm}
\ee
in order to model high-order multi-scale problems (see Chapter 4 \cite{gao-dual00} and \cite{gao-jogo06,gao-yu,gaot-ima}).

 In the set of 12 (= [6-14] + \cite{z16})  papers,  M.D. Voisei,
 C. \za and his former student  R. Strugariu  have made either mathematical mistakes (failed to correctly understand the canonical duality-triality theory and basic concepts in physics and nonlinear analysis),
 or ethic mistakes (repeatedly using the same
 condition $\dim \Pi \neq \dim \Pi^d$ as ``counter-examples"
 to against the double-min duality without citing author's original papers \cite{gao-opt03,gao-amma03},
 wherein this condition was first discovered).
 The mathematical mistakes show a huge gap between mathematical optimization and nonlinear analysis/mechanics.
 It is author's hope that by reading this paper, the readers can have  a clear  understanding
 not only on the canonical duality-triality theory and  its potential applications
  in multidisciplinary  fields,
but  also on the generalized  duality-triality principle and its role in modeling/understanding real-world problems. \\

\noindent{\bf Acknowledgements}:
The author should express his sincerely appreciation to
 his co-workers, students and colleagues for their understanding, excellent collaborations  and constant supports
for developing an unconventional theory in different fields, especially
  during the very difficult time caused by  the false challenges. 
Valuable comments and advices from  Professor Hanif Sherali at Virginia Tech and Professor Shu-Cherng Fang at North Carolina State University are highly appreciated.
 The author also would like to thank  Professor P. Krokhmal, the editor of {\em Optimization Letters},
 for his invitation to write  this paper.
This research has been continuously sponsored by US Air Force Office of Scientific Research under the grant AFOSR FA9550-10-1-0487 (2009-2015).
Strong supports from  program managers 
Dr. Jun Zhang, Dr. Jay Myung, and Dr. James Lawton are sincerely acknowledged.

\end{document}